\documentclass[11pt]{article}
\usepackage{amsmath,amsthm}
\usepackage{pstricks,pst-node,pst-tree}
\usepackage{latexsym}
\usepackage{amssymb,amscd}
\usepackage{graphics}   
\usepackage{url}
\usepackage[left=1.5in,top=1in,right=1in,nohead]{geometry}
\numberwithin{equation}{section}
\let \:=\colon
\let \beg=\begin
\let \mb=\mathbb
\let \mc= \mathcal
\let \ra=\rightarrow
\let \st=\stackrel
\let \sub=\subset
\let \Ga=\Gamma

\let \la=\lambda
\let \sig=\sigma
\let \al=\alpha
\let \be=\beta
\let \bet=\beta
\let \fl=\flushleft
\let \lfl=\lfloor
\let \rfl=\rfloor
\let \fr=\frac
\let \ga=\gamma

\let \om=\omega

\let \ov=\overline

\let \ub=\underbrace
\let \de=\delta
\let \De=\Delta
\let \Om=\Omega
\let \pr=\prime
\let \hra=\hookrightarrow
\let \lra=\longrightarrow

\let \ubr=\underbrace
\let \wt=\widetilde

\begin{document}
\title{\bf Geometry of curves with exceptional secant planes: linear series along the general curve}
\author{\small Ethan Cotterill}
\maketitle

\vspace{-20pt}
\beg{abstract}
We study linear series on a general curve of genus $g$, whose images are exceptional with
regard to their secant planes. Working in the framework of an extension of
Brill-Noether theory to pairs of linear series, we prove that a general
curve has no linear series with exceptional secant planes, in a very precise sense, whenever the total number of series is finite. Next, we partially solve the problem of computing the number of linear
series with exceptional secant planes in a one-parameter family in terms of
tautological classes associated with the family, by evaluating our hypothetical formula along judiciously-chosen test families. As an application, we compute the number of linear series with exceptional secant planes on a general curve equipped with a one-dimensional family of linear series. We pay special attention to the extremal case of $d$-secant $(d-2)$-planes to $(2d-1)$-dimensional series, which appears in the study of Hilbert schemes of points on surfaces. In that case, our formula may be rewritten in terms of hypergeometric series, which allows us both to prove that it is nonzero and to deduce its asymptotics in $d$.
\end{abstract}

\section{Introduction: Brill--Noether theory for \\pairs of series}
\vspace{-5pt}
Determining when an abstract curve $C$ comes equipped with a map to $\mb{P}^s$ of degree $m$ is central to curve theory. There is an {\it enumerative} aspect of this study, which involves determining formulas that describe the expected behavior of linear series along a curve. There is a {\it validative} aspect, which involves checking that the expected behavior holds. The Brill--Noether theorem, which is both enumerative and validative, asserts that when the Brill--Noether number $\rho(g,s,m)$ is nonnegative, $\rho$ gives the dimension of the space of series $g^s_m$ on a general curve $C$ of genus $g$, and that there is an explicit simple formula for the class of the space of linear series $G^s_m(C)$.

Since every linear series without base points determines a map to projective space, it is natural to identify a series with its image. Singularities of the image of a curve under the map defined by a series arise because the series admits certain subseries with base points; abusively, we refer to these subseries as ``singularities" of the series itself. Eisenbud and Harris \cite{EH1} showed that a general $g^3_m$ on a general curve of genus $g$ has no {\it double points}, or equivalently, that no inclusion
\vspace{-2pt}
\[
g^2_{m-2}+ p_1+ p_2 \hra g^3_m
\]
\vspace{-2pt}
exists, for any pair $(p_1,p_2)$ of points along the curve. They also showed that series with double points sweep out a divisor inside the space of all series $g^3_m$ along curves of genus $g$. 

In general, we say that an $s$-dimensional linear series $g^s_m$ {\it has a $d$-secant $(d-r-1)$-plane} provided an inclusion
\beg{equation}\label{seriesinseries}
g^{s-d+r}_{m-d}+ p_1 + \dots + p_d \hra g^s_m
\end{equation}
exists. The inclusion \eqref{seriesinseries} means that the image of the $g^s_m$ intersects a $(d-r-1)$-dimensional linear subspace of $\mb{P}^s$ in $d$-points; such a linear subspace is a ``$d$-secant $(d-r-1)$-plane". Hereafter, we use ``$d$-secant $(d-r-1)$-plane" to mean any inclusion \eqref{seriesinseries}. 

Next, let
\[
\mu(d,r,s):=d-r(s+1-d+r).
\]
The invariant $\mu$ computes the expected dimension of the space of $d$-secant $(d-r-1)$-planes along a fixed $g^s_m$. For example, when $\mu(d,r,s)=0$, we expect that the $g^s_m$ admits finitely many $d$-secant $(d-r-1)$-planes. 
Macdonald \cite{M} produced a general solution in the nineteen-fifties, though in practice his formulas are difficult to evaluate, as they require computing multi-indexed sums.

In this work, we study the analogous problem in case the $g^s_m$ is allowed to move, but the underlying curve is fixed. We show that a general curve admits {\it no} linear series with exceptional secant planes when $\rho=0$ and $\mu<0$. Similarly, we show that when $\rho=1$ and $\mu=-1$, there are finitely many series with exceptional secant planes along a general curve. 

We also develop a framework for computing the number of linear series with exceptional secant planes in a one-parameter family of series, based on a method of undetermined coefficients. In the present paper, we apply our method to produce formulas for the number of series with exceptional secant planes on a general curve when $\rho=1$ and $\mu=-1$. In a subsequent paper \cite{Co2}, we refine our technique to compute the classes of {\it secant plane} divisors on the moduli space $\ov{\mc{M}}_g$ associated to curves with linear series that are exceptional vis-\`{a}-vis their secant planes.

A couple of words are in order regarding the linear series parameter $r$. To avoid trivialities, we must have
\[
1 \leq r \leq s.
\]
Each specialization of $r$ defines an infinite family of examples, indexed by the incidence parameter $d$. For the applications to moduli, the most interesting aspect of our enumerative study concerns the large-$d$ asymptotics of our secant plane formulas, the behavior of which depends on the specialization we choose. The two most natural choices are $r=1$ and $r=s$, and in this paper we focus on the former, which corresponds to the situation studied by Lehn \cite{Le} in the context of the Hilbert scheme of points on a surface. Note that the case $r=1$ is ``tautological" in that it corresponds to the situation in which the evaluation map
\[
V \ra H^0(L/L(-p_1-\dots-p_d))
\]
corresponding to a given linear series $(L,V)$ fails to be surjective along a $d$-tuple of points $p_1, \dots, p_d$ on the curve in question.
We show that when $r=1$, our formulas can be compactly expressed in terms of generating functions, each term of which is a finite linear combination of hypergeometric functions of type ${\empty}_3F_2$.

\section*{Acknowledgements}
This work is part of my doctoral thesis, which was carried out under the supervision of Joe Harris. I thank Joe, and Steve Kleiman at MIT, for many valuable conversations. I also thank the referee for a very careful reading of the paper, and valuable comments which led to a substantial improvement in the exposition.

\subsection{Roadmap}
The material following this introduction is organized as follows. In the second section, we address the validative problem of determining when a curve possesses linear series with exceptional secant planes.  Our first two results establish that on a general curve, there are no linear series with exceptional secant planes when the expected number of such series is zero. We begin by proving the following basic nonexistence result:

{\fl \bf Theorem 1.}
{\it If $\rho=0$ and $\mu$ is negative, then a general curve $C$ of
genus $g$ admits no $s$-dimensional linear series $g^s_m$ with
$d$-secant $(d-r-1)$-planes.}

To prove Theorem 1, we show that on a particular ``nearly"-stable model of a $g$-cuspidal rational curve (obtained by blowing up the stable model in finitely many points), there are no linear series with exceptional secant planes when $\rho=0$ and $\mu=-1$. In \cite{Fa2}, which appeared as a preprint at the same time that an announcement of the results in the current paper was circulating, G. Farkas obtains a proof of the natural generalization of Theorem 1 to the case $\rho \geq 0$ via limit linear series. The argument which we present is substantially simpler, if less far-reaching, and naturally generalizes the argument used in \cite[Prop. 5.52]{HM} to show that no linear series exist on a general curve when $\rho$ is negative. Our argument also highlights the r\^{o}le of the two-step flag variety in these questions.

Finally, we close the second section by proving the following theorem, which gives geometric significance to the enumerative study carried out in Section 4:
{\fl \bf Theorem 2.}
{\it If $\rho=1$ and $\mu=-1$, then there are finitely many linear series $g^s_m$ with $d$-secant $(d-r-1)$-planes on a general curve $C$ of genus $g$.}

\medskip
In Section 3, we begin our enumerative study of linear series with exceptional secant planes along a general curve $C$. We start by considering the more general problem of computing the expected number of linear series with exceptional secant planes in a one-parameter family of linear series (not necessarily along a fixed curve) by computing the number of exceptional series along judiciously-chosen ``test families". In general, we know from Ran's work \cite{R2} that for a one-parameter family $\pi: \mc{X} \ra B$ of linear series $g^s_m$, the number of fibers $N^{d-r-1}_d$ with $d$-secant $(d-r-1)$-planes is given by a universal formula
\beg{equation}\label{univpoly}
N^{d-r-1}_d= P_{\al} \al+ P_{\be} \be+ P_{\ga} \ga+ P_c c + P_{\de_0} \de_0
\end{equation}
where $\al, \be, \ga, c$, and $\de$ are tautological numerical invariants associated to $\pi$, and the coefficients $P_{\al}, P_{\be}$, $P_{\ga}, P_c$, and $P_{\de_0}$ are polynomials in $d, m, r$, and $s$. Here \eqref{univpoly} holds whenever the number of such fibers is finite. However, whenever $s \geq 3$ and $C$ is general in moduli, we have $\ga= \de_0=0$, so only three tautological relations are required (instead of five, in the more general setting). Section 3.3 is devoted to establishing the enumerative nature of our two most basic relations among tautological coefficients, which are derived from the study of the enumerative geometry of a {\it fixed} curve in projective space carried out in \cite{ACGH}.

When $r=1$, our results are strongest. A key ingredient in these is the generating function for the expected number $N_d$ of $d$-secant $(d-2)$-planes to a $g^{2d-2}_m$, which we obtain in  Section~\ref{Ndfmla}. We show:
{\fl \bf Theorem~\ref{Ndthm}.}
\[ 
\sum_{d \geq 0} N_d z^d= \biggl(\fr{2}{(1+4z)^{1/2}+1} \biggr)^{2g-2-m} \cdot (1+4z)^{\fr{g-1}{2}}.
\]
Lehn \cite{Le} studied the analogous problem for Hilbert schemes of points on surfaces using representation-theoretic methods, and obtained a conjectural generating function for the corresponding numbers $N_d$. As is clear from the proof of Theorem~\ref{Ndthm}, $d$-secant $(d-2)$-planes to linear series on curves are closely related to Catalan numbers, whose generating series $C(z)$ satifies $C(-z)=\fr{2}{(1+4z)^{1/2}+1}$. While revising this paper, the author learned that Le Barz \cite{Lb2} has also recently proved Theorem~\ref{Ndthm}, as a consequence of the algorithm for computing the multisecant loci of a fixed curve described in \cite{Lb1}. Our methods are more elementary than Le Barz's. The combinatorics of $d$-secant $(d-2)$-planes explored in this paper and in the thesis \cite{Co1} has led to interesting new combinatorial identities \cite{DY, SZ}.

In Section~\ref{genfunctions}, we use the generating function for $N_d$ obtained in Section~\ref{Ndfmla} to determine generating functions for the tautological coefficients $P$, whenever $r=1$. In Section~\ref{hypergeomfns}, we use the generating functions determined in Section~\ref{genfunctions} in order to realize (in Theorem 4) each of the tautological coefficients $P_{\al}, P_{\be}, \text{ and } P_c$ as linear combinations of generalized hypergeometric functions.

Finally, in Section~\ref{rho=1} we determine an enumerative formula for the number of linear series with exceptional secant planes along a general curve when $\rho=1$. Namely, we prove:
{\fl \bf Theorem~\ref{genericcurve}.}
{\it
Let $\rho=1, \mu=-1$. The number $N^{\pr, d-r-1}_d$ of linear series $g^s_m$ with $d$-secant $(d-r-1)$-planes on a general curve of genus $g$ is given by
\vspace{-3pt}
\[
\beg{split}
N^{\prime, d-r-1}_d= &\fr{(g-1)! 1! \cdots s!}{(g-m+s)! \cdots (g-m+2s-1)!(g-m+2s+1)!} \cdot \\
&[(-g m + m^2  - 3 m s + 2 s^2  - m + s - g)A\\
&+(g d + g - m d - m + 2 s d + 2 s + d + 1)A^{\prime}]
\end{split}
\]
}
where $A$ and $A^{\pr}$ compute, respectively, the expected number of $d$-secant $(d-r)$-planes to a $g^{s+1}_m$ that intersect a general line, and the expected number of $(d+1)$-secant $(d-r)$-planes to a $g^{s+1}_{m+1}$. Note that formulas for $A$ and $A^{\pr}$ were computed by Macdonald in \cite{M}.

Subsequently, we specialize to the case $r=1$, where we obtain a hypergeometric formula for the number $N^{\pr, d-2}_d$ of $(2d-1)$-dimensional series with $d$-secant $(d-2)$-planes along a general curve when $\rho=1$. Using that formula, we prove Theorem~\ref{ndprimepos}, which characterizes exactly when $N^{\pr, d-2}_d$ is positive, and we determine the asymptotics of $N^{\pr, d-2}_d$ as $d$ approaches infinity.

\section{Validative study}\label{validative}
We begin by proving the following theorem.
\beg{thm}\label{thm1}
Assume that $\rho=0$, and $\mu$ is negative. Under these conditions, a general curve $C$ of
genus $g$ admits no $s$-dimensional linear series $g^s_m$ with
$d$-secant $(d-r-1)$-planes.
\end{thm}
The theorem asserts that on $C$, there are no {\it pairs} of series
$(g^{s-d+r}_m,g^s_m) \in G^{s-d+r}_m(C) \times G^s_m(C)$ satisfying \eqref{seriesinseries} for any choice of $d$ base points $(p_1,\dots,p_d) \in C^d$.
To prove it, we specialize $C$ to a broken flag curve $\widetilde{C}$ of the type
used in Eisenbud and Harris' proof of the Giesker-Petri
theorem \cite{EH2}: $\widetilde{C}$ is a nodal curve of compact type comprised of a ``spine" of rational curves
$Y_i$, some of which are linked via a sequence of rational curves to $g$ elliptic ``tails'' $E_1,\dots,E_g$.  See Figure 1. The components $Y_i$ of the spine are numbered so that the index $i$ increases as one traverses the spine from the top of the figure to the bottom. We set $q_i:= Y_{i-1} \cap Y_i$ for all $i \geq 2$.

It suffices to show that $\widetilde{C}$ admits no inclusions \eqref{seriesinseries} of {\it limit linear series}. So assume for the sake of argument that $\widetilde{C}$ {\it does} in fact admit a (limit linear) series $g^{s-d+r}_m
\hra g^s_m$ satisfying \eqref{seriesinseries}. We will obtain a
contradiction by showing that \eqref{seriesinseries} is incompatible with
basic numerical restrictions obeyed by the vanishing sequences of the
$g^s_m$ and $g^{s-d+r}_m$ at intersection points of rational components along the
spine.

In what follows, let $V_Z$ denote the aspect of the $g^s_m$ along the component $Z
\subset \widetilde{C}$. We will systematically use the following three basic
facts from the theory of limit linear series \cite{EH3}:
\beg{itemize}
\item {\bf LS1.} At a node $q$ along which components $Y, Z \subset \widetilde{C}$
  intersect transversely, the vanishing sequences $a(V_Y,q)$ and $a(V_Z,q)$ verify
\[
a_i(V_Y,q)+a_{s-i}(V_Z,q) \geq m
\]
for all $0 \leq i \leq s$.
\item {\bf LS2.} Assume that a set of compatible
  bases for $V$ along $\widetilde{C}$ has been chosen, in the sense that
  $V_{Y_i} \subset V_{Y_{i+1}}$, for every $i$. Then
 \[  
  a_j(V_{Y_{i+1}},q_{i+1}) \geq a_j(V_{Y_i},q_i).
\]
for every index $0 \leq j \leq s$.
\item {\bf LS$2^{\pr}$.} If $\rho(g,s,m)=0$, the following strengthening of LS2 holds.
  \beg{itemize} 
  \item If $Y_i$ is linked via rational curves to an elliptic tail, then
  \[
  a_j(V_{Y_{i+1}},q_{i+1})=a_j(V_{Y_i},q_i)+1
  \] 
  for all $0 \leq j \leq s$ except for
  a single index $j$, for which 
  \[
  a_j(V_{Y_i},q_{i+1})=a_j(V_{Y_i},q_i).
  \] 
  \item If $Y_i$ is not linked via rational curves to an elliptic tail, then 
\[  
  a(V_{Y_{i+1}},q_{i+1})=a(V_{Y_i},q_i).
\]
  \end{itemize}

 \item {\bf LS3.} The ramification sequence of a $g^s_m$ at a point $p \in \mb{P}^1$ determines a Schubert variety in $\mb{G}(s,m)$. Intersections of (arbitrarily many) Schubert varieties associated to ramification sequences $\alpha_i=\alpha(V,r_i)$ at distinct points $r_i$ have the expected dimension. Consequently, a smooth rational curve admits a $g^s_m$ with ramification
  sequences $\alpha_i$ at $r_i$
  if and only if the product of the corresponding Schubert cycles
is nonzero in $H^*(\mb{G}(s,m),\mb{Z})$. The $r_i$ need not be general in order for ``dimensional transversality" to hold \cite[Thm. 2.3]{EH1}.
\item {\bf LS4.} Let $(L,V)$ denote a linear series along a reducible curve $Y \cup_q Z$. If $Z$ is a smooth and irreducible elliptic curve, then the aspect $V_Z$ of the linear series along $Z$ has a cusp at $q$, i.e., the ramification sequence $\al(V_Z,q)$ satisfies
\[
\al(V_Z,q) \geq (0,1,\dots,1).
\]
\end{itemize}

\beg{figure}\label{Fig1}
\beg{picture}(40,180)
\qbezier(0,180)(20,160)(0,145) 
\qbezier(0,150)(20,130)(0,110) 
\put(0,130){$Y_i$}
\put(-5,95){$Y_{i+1}$}
\qbezier(0,125)(20,115)(40,125)
\qbezier(30,125)(50,115)(70,125)
\qbezier(60,125)(80,115)(100,125)
\qbezier(90,125)(110,115)(130,125)
\qbezier(120,125)(180,115)(340,115)
\qbezier(0,120)(20,100)(0,80) 
\qbezier(0,70)(20,60)(40,70)
\qbezier(30,70)(50,60)(70,70)
\qbezier(60,70)(80,60)(100,70)
\qbezier(90,70)(110,60)(130,70)
\qbezier(120,70)(180,60)(255,60)
\qbezier(0,60)(20,40)(0,20) 
\qbezier(0,30)(20,10)(0,-10) 
\qbezier(0,00)(20,-20)(0,-40) 
\qbezier(0,10)(20,0)(40,10)
\qbezier(30,10)(50,0)(70,10)
\qbezier(60,10)(80,0)(100,10)
\qbezier(90,10)(110,0)(130,10)
\qbezier(120,10)(180,0)(340,0)
\put(30,90){\vector(-1,0){20}}
\put(30,90){\vector(-1,-4){18}}
\put(30,90){\vector(-1,4){18}}
\put(33,98){Spine of rational curves}
\put(270,60){\vector(-1,4){14.5}}
\put(270,60){\vector(-1,0){12.5}}
\put(270,60){\vector(-1,-4){14.5}}
\put(270,70){Elliptic}
\put(272.5,57.5){tails}
\put(-5,75){\dots}
\put(-5,63){\dots}
\put(270,119){$E_1$}
\put(270,5){$E_g$}

\qbezier(292.5,60)(311,60)(340,60)
\end{picture}
\caption{A broken flag curve.}
\end{figure}

By repeated blowing-up, we are free to assume that no base point $p_i$ lies at a point of attachment linking components of $\widetilde{C}$.
For convenience, we also make the following simplifying assumption, which we will remove later:

{\it No $p_i$ lies along an elliptic tail.}

Now fix a component $Y_i$ along the spine. If it is interior to the spine, then it has two ``special" points corresponding to the intersections with adjacent rational components $Y_{i-1}$ and $Y_{i+1}$ along the spine, which we mark by $0$ and $\infty$, respectively. If it is linked via a chain of rational curves to an elliptic tail, then we denote its intersection with the first component of the chain by $1$. If $Y_i$ is not interior to the spine, and is marked by $0$ (resp., $\infty$), then we mark an additional point by $\infty$ (resp., $0$), so that every component along the spine of $Y$ has at least two marked points. 

Denote the {\it vanishing orders} of $V_{Y_i}$ at $0$ (resp., $\infty$) by $a_j$ (resp., $b_j$), $0 \leq j \leq s$; if $V_{Y_i}$ is spanned by sections $\sig_j(t), 0 \leq j \leq s$ in a local uniformizing parameter $t$ for which $\mbox{ord}_t(\sig_i)< \mbox{ord}_t(\sig_j)$ whenever $i<j$, then $a_i:=\mbox{ord}_t(\sig_i)$. Denote the corresponding vanishing orders of the $g^{s-d+r}_m$ along $Y_i$ by $u_j$ and $v_j$, respectively. Note that the vanishing sequence $(u_j)$ (resp., $(v_j)$) is a subsequence of $(a_j)$ (resp, $(b_j)$). Recall that $(u_j)$ and $(v_j)$ correspond to Schubert cycles $\sig_{\al}$ and $\sig_{\be}$ in $H^*(\mb{G}(s-d+r,m),\mb{Z})$ defined by
\[
\al_j= u_{s-d+r-j}-(s-d+r)+j, \text{ and } \be_j= v_{s-d+r-j}-(s-d+r)+j
\]
for all $j=0, \dots, s-d+r$, respectively. 


We say that the vanishing sequences $(u_j)$ and $(v_j)$ are {\it complementary} if 
\[
u_j=a_{k(j)} \text{ and } v_j=b_{s-k(s-d+r-j)}
\]
for some sequence of nonnegative integers $k(j), j=0, \dots, s-d+r$. Base points cause $(u_j)$ and $(v_j)$ to fail to be complementary to one another by an amount which may be estimated uniformly.

\beg{lem}\label{lem1}
Assume that the $g^{s-d+r}_m$ along $Y_i$ has a simple base base point $p$. Then
\[
v_j=b_{s-k(s-d+r-j)-k^{\prime}(j)}, j=0, \dots, s-d+r
\]
for some sequence of nonnegative integers $k^{\prime}(j), j=0,\dots,s-d+r$, at least $(s-d+r)$ of which are equal to at least 1.
\end{lem}

In other words, a simple base point leads to $(s-d+r)$ ``shifts" of vanishing order indices of our $g^{s-d+r}_m$.

{\fl \bf Remark.} Inclusions of linear series $g^{s-d+r}_m \hra g^s_m$ on $\mb{P}^1$ determine a two-step flag variety $\mbox{Fl}(s-d+r,s; m)$, which comes equipped with a natural projection to $\mb{G}(s,m)$, the fibers of which are isomorphic to $\mb{G}(s-d+r,s)$. A pair of vanishing sequences for $g^s_m$ and its included $g^{s-d+r}_m$ at a point defines a Schubert variety in $\mbox{Fl}(s-d+r,s; m)$. The codimension of this Schubert variety is the sum of two terms:
\beg{enumerate}
\item[(i)] the codimension of its image in $\mb{G}(s,m)$, i.e., $\sum_{j=0}^s (a_j-j)$, where $(a_0,\dots,a_s)$ is the vanishing sequence of the $g^s_m$;
\item[(ii)] the total number of shifts $\sum_{j=0}^{s-d+r} (k(j)-j)$, where $(a_{k(0)},\dots, a_{k(s-d+r)})$ is the vanishing sequence of the included $g^{s-d+r}_m$.
\end{enumerate}
Complementarity encodes the notion of ``maximal codimension in a fiber" over a point of $\mb{G}(s,m)$. That is, for a fixed choice of series $g^s_m$ together with the choice of vanishing sequences at 0 and $\infty$, we expect there to be finitely many subseries $g^{s-d+r}_m$ with complementary vanishing subsequences at 0 and $\infty$.

{\fl \bf Proof of Lemma 1.}
Set $Y:=Y_i$. To increase readability, we use $\sig(x_0,\dots,x_{s-d+r})$ to denote the Schubert cycle $\sig_{x_0,\dots,x_{s-d+r}}$.

\medskip
{\it Subcase: $Y$ is linked, via a chain of rational curves, to an elliptic tail.}

Recall that we denote the point of attachment of the chain along $Y$ by 1. By LS4, the aspect along $Y$ of the $g^{s-d+r}_m$ has at least a cusp at $1$; i.e., the corresponding Schubert cycle $\sig^{(1)}$ satisfies
\[
\sig^{(1)} \geq \sig_{1,\dots,1,0}.
\]
Meanwhile, the Schubert cycle corresponding to $p$ in $H^*(\mb{G}(s-d+r,m),\mb{Z})$ is $\sig^{(p)}= \sig_{1,\dots,1}$.
{\fl By} LS3, the intersection
\beg{equation}\label{schub1}
\sig^{(0)} \cdot \sig^{(1)} \cdot \sig^{(\infty)} \cdot \sig^{(p)} \in H^*(\mb{G}(s-d+r,m),\mb{Z})
\end{equation}
is necessarily nonzero.
Since $\sig^{(p)}= \sig_{1,\dots,1}$, \eqref{schub1} is clearly nonzero if and only if the corresponding intersection
\[
\sig^{(0)} \cdot \sig^{(1)} \cdot \sig^{(\infty)}
\]
is nonzero in $H^*(\mb{G}(s-d+r,m-1),\mb{Z})$. In particular, we must have
\[
\sig^{(0)} \cdot \sig^{(\infty)} \cdot \sig_{1,\dots,1,0} \neq 0 \in H^*(\mb{G}(s-d+r,m-1),\mb{Z}).
\]

{\fl Now say} that the vanishing sequence of the $g^s_m$ along $Y$ at $0$ is 
\[
a(V_Y,0)=(a_0,\dots,a_s)
\]
and that, correspondingly, the vanishing sequence of the $g^{s-d+r}_m$ at $0$ is
\[
(u_0,\dots, u_{s-d+r})=(a_{k(0)},\dots,a_{k(s-d+r)})
\]
for some sequence of nonnegative integers $k(j), j=0, \dots, s-d+r$.
We then have
\[
\sig^{(0)}= \sig(a_{k(s-d+r)}-(s-d+r),\dots,a_{k(1)}-1,a_{k(0)}).
\]
The sequence
\[
(b_{s-k(s-d+r)},\dots,b_{s-k(0)})
\]
is complementary to $(u_0,\dots, u_{s-d+r})$. Let $\sig^{(0^{\vee})}$ denote the corresponding Schubert cycle; then
\[
\sig^{(0^{\vee})}= \sig(b_{s-k(0)}-(s-d+r),\dots,b_{s-k(s-d+r)}).
\]

{\fl \bf A key observation.} Combining LS1 with LS$2^{\pr}$, we obtain
\beg{equation}\label{keyobs}
b_{s-i}=m-1-a_i
\end{equation}
for every $i$ in $\{0,\dots,s\}$, except for a unique index $j$ for which $b_{s-j}=m-a_j$. It follows immediately that the intersection product
\[
\sig^{(0)} \cdot \sig^{(0^{\vee})} \in H^*(\mb{G}(s-d+r,m-1))
\]
is zero unless 
\[
b_{s-k(j)}=m-1-a_{k(j)}
\]
for all $j$ in $\{0,\dots,s-d+r\}$.

{\fl On} the other hand, by Pieri's rule
\[
\sig^{(0)} \cdot \sig_{1,\dots,1,0}
\]
is a sum of Schubert cycles 
\[
\sig^{(0^{\pr})}= \sig(a_{k(s-d+r)}-(s-d+r)+ k^{\pr \pr \pr }(s-d+r),\dots,a_{k(1)}-1+ k^{\pr \pr \pr}(1),a_{k(0)}+ k^{\pr \pr \pr}(0))
\]
for some sequence of nonnegative integers $k^{\pr \pr \pr }(j), 0 \leq j \leq s-d+r$, at least $(s-d+r)$ of which are equal to at least one.

{\fl Now} say that the vanishing sequence of the $g^{s-d+r}_m$ at $\infty$ is 
\[
a(V_Y, \infty)= (b_{s-k(s-d+r)}- k^{\pr}(0), \dots, b_{s-k(0)}- k^{\pr}(s-d+r))
\]
for some sequence of nonnegative integers $k^{\pr}(j), j=0, \dots, s-d+r$.

The corresponding Schubert cycle is
\[
\sig^{(\infty)}=\sig(b_{s-k(0)}-(s-d+r)-k^{\pr}(0),\dots,b_{s-k(s-d+r)}-k^{\pr}(s-d+r)).
\]

If the intersection product
\[
\beg{split}
\sig^{(0^{\pr})} \cdot \sig^{(\infty)}&= \sig(a_{k(s-d+r)}-(s-d+r)+ k^{\pr \pr \pr }(s-d+r),\dots,a_{k(0)}+ k^{\pr \pr \pr}(0)) \\ &\cdot \sig(b_{s-k(0)}-(s-d+r)-k^{\pr}(0),\dots,b_{s-k(s-d+r)}-k^{\pr}(s-d+r))
\end{split}
\]
is nonzero, then the $(s-d+r+1)$ sums of complementary indices
\[
\beg{split}
&a_{k(s-d+r)}-(s-d+r)+ k^{\pr \pr \pr }(s-d+r)+ b_{s-k(s-d+r)}-k^{\pr}(s-d+r) \\
& \dots \\
&a_{k(0)}+ k^{\pr \pr \pr}(0)+ b_{s-k(0)}-(s-d+r)-k^{\pr}(0)
\end{split}
\]
are each at most $m-1-(s-d+r)$. Via \eqref{keyobs}, we conclude that
\[
\beg{split}
m-1 - (s-d+r)+ k^{\pr \pr \pr}(j)- k^{\pr}(s-d+r-j) &\leq m-1-(s-d+r), \text{ i.e., that}\\
k^{\pr \pr \pr}(j) \leq k^{\pr}(s-d+r-j)
\end{split}
\]
for every $j$ in $\{0, \dots, s-d+r\}$.

{\fl But} $(s-d+r)$ values of $k^{\pr \pr \pr}(j)$ are nonzero, from which it follows that the same is necessarily true of the values of $k^{\pr}(j)$. The desired conclusion follows immediately.

{\it Subcase: $Y$ is not linked to an elliptic tail.} 

In this case, by \cite[Lem 5.57, pt 1]{HM}, the aspect of the $g^s_m$ along $Y=Y_i$ satisfies
 \[
  a_j(V_{Y_{i+1}},q_{i+1})=a_j(V_{Y_i},q_i)+1
  \] 
  for all $0 \leq j \leq s$. Accordingly, LS2 implies that
\beg{equation}\label{keyobsbis}
b_{s-i}= m- a_i
\end{equation}
for every $i$ in $\{0,\dots,s\}$. By arguing as in the first subcase using \eqref{keyobsbis} in place of \eqref{keyobs}, we deduce that $k^{\pr}(j) \geq 1$ for every $j$ in $\{0, \dots, s-d+r\}$. $\Box$




{\fl \bf Proof of Theorem~\ref{thm1}.}

\medskip
{\it Proof in case all base points of the limit series $g^{s-d+r}_m$ are simple.}

Let $Y_i=Y$ denote the same component of $\wt{C}$ as before, with marked points 0 and $\infty$. Recall that $(a_j)$ (resp., $(u_j)$) denotes the vanishing sequence of the aspect of the $g^s_m$ (resp., $g^{s-d+r}_m$) at 0. Similarly, $(b_j)$ (resp., $(v_j)$) denotes the vanishing sequence of the aspect of the $g^s_m$ (resp., $g^{s-d+r}_m$) at $\infty$. Now define the sequence $(u_j^{\prime}), j=0, \dots, s-d+r$ by setting $u_j^{\prime}:= m- v_j$ for every $j$. If $Y_i$ is interior to the spine of $\widetilde{C}$, then, by LS1, $(u_j^{\prime})$ is a subsequence of the vanishing sequence $a(V_{Y_{i+1},p_{i+1}})=(a^{\pr}_0,\dots,a^{\pr}_s)$. Letting
\[
u_j^{\pr}= a^{\pr}_{k^{\pr \pr}(j)},
\]
the first lemma asserts that the sequences $k(j)$ and $k^{\pr \pr}(j)$ satisfy
\[
k^{\pr \pr}(j) \geq k(j)+1
\]
for at least $(s-d+r)$ values of $j$. In other words, in passing from $Y_i$ to $Y_{i+1}$ the existence of a simple base point for the $g^{s-d+r}_m$ forces $(s-d+r)$ vanishing order indices $k(j)$ to ``shift to the right" by at least one place. Similarly, the existence of $d$ simple base points along distinct components of $\wt{C}$ forces $d(s-d+r)$ shifts. On the other hand, shifts of vanishing order indices are constrained; namely, each index can shift at most $s-(s-d+r)=d-r$ places. So the maximum possible number of shifts is $(s-d+r+1)(d-r)$, and we necessarily have
\[
d(s-d+r) \leq (s-d+r+1)(d-r),
\]
which contradicts the fact that
$\mu(d,r,s)=(s-d+r+1)(d-r)- d(s-d+r)= -1$.

Thus, we conclude whenever all base points occur along interior components of the spine of $\widetilde{C}$. An analogous argument yields a contradiction whenever (simple) base points lie along either of the two ends of the spine. Namely, assume that $Y$ is the ``bottommost" component of the spine, so that $Y$ is marked by $0$ (the intersection with the component of the spine directly above it) and $\infty$. Then $(d-1)$ base points of the $g^{s-d+r}_m$ lie on components above $Y$, forcing $(d-1)(s-d+r)$ shifts of vanishing order indices. The vanishing sequence of the $g^{s-d+r}_m$ at 0, viewed as a subsequence of $(a_j)=a(V_{Y},0)$, is
\[
(u_j)= (a_{k(j)}), j=0, \dots, s-d+r
\]
where $\sum_{j=0}^{s-d+r} k(j)-j \geq (d-1)(s-d+r)$. Applying Lemma 1, we deduce that the vanishing sequence $(v_j)=(b_{s-k(s-d+r-j)-k^{\pr}(j)})$, viewed as subsequence of $(b_j)=a(V_{Y},\infty)$, satisfies
\[
k^{\pr}(j) \geq 1
\]
for at least $(s-d+r)$ indices $j \in \{0, \dots, s-d+r \}$. But then
\[
\sum_{j=0}^{s-d+r} [k^{\pr}(j)+ k(s-d+r-j)- (s-d+r-j)] \geq d(s-d+r),
\]
i.e., $d(s-d+r)$ shifts of vanishing order indices are forced at $\infty$, which is impossible.

\medskip
{\it Proof of Theorem~\ref{thm1}, assuming all base points lie along the spine.}

By blowing up, if necessary, we may assume that there is a single base point $p_i$ (possibly multiple) along each component. 

Now fix a component $Y$ of $\wt{C}$ along which the $g^{s-d+r}_m$ has a multiple base point $d^{\pr}p, d^{\pr} >1$. Let 0 and $\infty$ designate the same marked points along $Y$ as before, with $((u_j),(a_j))$ and $((v_j),(b_j))$ the corresponding vanishing sequences of the $g^{s-d+r}_m$ and $g^s_m$. We further assume that 1 is marked, i.e., that $V_Y$ has a cusp there (the proof is easier when 1 is unmarked.) As before, we have
\[
u_j=a_{k(j)}, \text{ and } v_j= b_{s-k(s-d+r-j)-k^{\pr}(j)}
\]
for certain sequences of nonnegative integers $k(j)$ and $k^{\pr}(j)$, $j=0, \dots, s-d+r$. 

\medskip
We need the following generalization of Lemma 1.
\beg{lem}\label{lem2}
Assume that the $g^{s-d+r}_m$ along $Y_i$ has a base point $d^{\pr}p$, $d^{\pr} \geq 1$. Then
\[
v_j=b_{s-k(s-d+r-j)-k^{\prime}(j)}, j=0, \dots, s-d+r
\]
for some sequence of nonnegative integers $k^{\prime}(j), j=0,\dots,s-d+r$, at least $(s-d+r)$ of which are equal to at least 1.
\end{lem}

{\bf \fl Proof of Lemma 2.} Set $Y=Y_i$, as before.

{\it Subcase: $s-d+r=1$, i.e., $g^{s-d+r}_m$ is a subpencil of $g^s_m$.} 

It suffices to show that $d^{\pr}$ of shifts of vanishing order indices are forced by $d^{\pr} p$, i.e., that 
\[
\sum_{j=0}^1 k^{\pr}(j) \geq d^{\pr};
\]
the proof of Theorem 1 then proceeds as in the case of simple base points. 

For this purpose, begin by fixing a basis $\{\sig_0, \dots, \sig_s\}$ for the $s$-dimensional series $V_Y$, such that the vanishing order of $\sig_j$ at 0 is $a_j(V_Y,0)$, for $j=0, \dots, s$. A convenient choice, given in \cite[p. 74]{EH4}, is
\beg{equation}\label{rhozeropres}
\sig_i= x_i t^{a_i} u^{m-a_i}+ t^{a_i+1}u^{m-a_i-1} \text{ for all }i \neq j, \sig_j= t^{a_j}u^{m-a_j}. 
\end{equation}
Here $(t,u)$ are homogeneous coordinates on $\mb{P}^1$ with respect to which $0= (0,1)$, $\infty= (1,0)$, and $1=(1,1)$, and $x_i=(1-a_j+a_i)/(a_j-a_i)$.
 
Now form a set $\wt{\mc{B}}_Y$ in the following way. To $\{\sig_{k(0)}, \sig_{k(1)}\}$, add sections $\sig_{k(j)+1}, \dots, \sig_{k(j)+k^{\pr}(j)}$ for each $j \in \{0,1\}$ with $k^{\pr}(j)>0$. Let $N$ denote the number of shifts of vanishing order indices which are forced by $d^{\pr} p$. The corresponding $(N+1)$-dimensional subseries $\wt{V}_Y \sub V_Y$ contains the $g^1_m$ of interest to us. 

On the other hand, $\wt{V}_Y$ has at least a cusp at $1$, and each of the sections in $\wt{\mc{B}}_Y$ vanishes to total order at least $(m-1)$ at $\{0 \cup \infty\}$. Whence, the total ramification degree of $\wt{V}_Y$ away from $\{0,1,\infty\}$ is at most 1.

Note that degree-$m$ subpencils of $\wt{V}_Y$ determine a $\mb{G}(1,N+1)$. If $V_Y$ is unramified at $p$, then pencils with a $d^{\pr}$-fold base point comprise a subvariety $\mc{W}$ of codimension $2d^{\pr}$. If $V_Y$ is simply ramified at $p$, $\mc{W}$ has codimension $2d^{\pr}-1$. As $\mc{W}$ is nonempty by assumption and $\mb{G}(1,N+1)$ is $2N$-dimensional, this forces $N \geq d^{\pr}$.

\medskip
{\it General case.} Each subpencil of the $g^{s-d+r}_m$ along $Y$ under consideration has a $d^{\pr}$-fold base point at $p$. Consequently, such a subpencil induces $d^{\pr}$ shifts of vanishing order indices. By varying the choice of subpencil, we conclude that the $d^{\pr}$-fold base point of the $g^{s-d+r}_m$ along $Y$ forces $d^{\pr}(s-d+r)$ shifts of vanishing order indices, and the proof of Theorem 1 now proceeds as in the case of simple base points.

\medskip
{\bf \fl Conclusion of the proof of Theorem 1}.

Finally, we explain how to remove the simplifying assumption inserted at the beginning. Namely, assume that the $g^{s-d+r}_m$ admits a base point $d^{\pr}p$ along an elliptic tail $E$, where $d^{\pr} \geq 1$. Say that $E$ intersects the rational component $Z$ of $\widetilde{C}$ in a node $q$ of $\widetilde{C}$. Note that the vanishing sequence at $q$ of the $g^{s-d+r}_m$ along $E$ is bounded above by
\[
(m-s+d-r-1-d^{\pr}, \dots, m-2-d^{\pr},m-d^{\pr});
\]
otherwise, the subpencil of sections of the $g^{s-d+r}_m$ along $E$ that vanish to maximal order define (upon removal of the $(m-d^{\pr}-2)$-fold base point $(m-d^{\pr}-2)q$) a $g^1_1$, which is absurd. It follows, by LS1, that the vanishing sequence at $q$ of the $g^{s-d+r}_m$ along $Z$ is {\it at least}
\[
(d^{\pr},d^{\pr}+2,\dots, s-d+r+1+d^{\pr}),
\]
which in turn implies that the same estimate holds for the vanishing sequence of the $g^{s-d+r}_m$ along the rational component $Y_i$ of the spine of $\widetilde{C}$ linked to $E$ at the corresponding node $\widetilde{q}$.

In other words, if the $g^{s-d+r}_m$ has a base point along $E$, then the $g^{s-d+r}_m$ also has a base point and a cusp along $Y_i$. In this way, we are reduced to the simplified setting in which no base points $p$ of the included series $g^{s-d+r}_m$ lie along elliptic components of $\widetilde{C}$, and are free to argue as before. \qed

{\bf \fl Remarks.}
\beg{itemize}
\item Let $Y \cong \mb{P}^1$ be a smooth rational curve with marked points $0, 1, \text{ and }\infty$. Let \\$\rho_Y(\al; \be; \ga; s,m)$ denote the dimension of the space of series $g^s_m$ along $Y$ with ramification sequence $(\al_j= a_j-j)$ (resp., $(\be_j=b_j-j), \ga$) at 0 (resp., $\infty, 1$). Let $V_Y$ be any $g^s_m$ with ramification $\al$, $\be$, and $(0,1,\dots,1)$ at 0, 1, and $\infty$; assume that $V_Y$ admits an inclusion
\[
g^{s-d+r}_{m-d^{\pr}}+ d^{\pr}p \hra g^s_m
\]
with $d^{\pr}\geq 1$ and $p \notin \{0,1,\infty\}$.
In the general case $\rho \geq 0$, one expects the following result to hold.

{\it
The vanishing sequence $(b_{s-k(s-d+r-j)-k^{\pr}(j)})_{j=0}^{s-d+r}$ of the included $g^{s-d+r}_m$ at $\infty$ is shifted at least $d^{\pr}(s-d+r)- \rho_Y(\al; \be; (0,1,\dots,1); s,m)$ places relative to the sequence complementary to $a(g^{s-d+r}_m,0)=(a_{k(j)})_{j=0}^{s-d+r}$. That is,
\beg{equation}\label{basicest}
\sum_{k=0}^{s-d+r} k^{\pr}(j) \geq d^{\pr}(s-d+r)- \rho_Y(\al; \be; (0,1,\dots,1); s,m).
\end{equation}
}

\item Now let $\wt{C}$ be a flag curve equipped with a $g^s_m$. By additivity of the Brill--Noether number, $\rho=\rho(g,s,m)$ is equal to $\rho_{Y_i}(\al; \be; (0,1,\dots,1); s,m)$ over all spinal components $Y_i$ of $\wt{C}$. In particular, the basic estimate \eqref{basicest} generalizes Lemma 2, as the number $\rho_Y(\al; \be; (0,1,\dots,1); s,m)$ is zero when $\rho=0$. 

To prove that $\wt{C}$ (and therefore, a general curve of genus $g$) admits no inclusion \eqref{seriesinseries} when $\rho + \mu < 0$ for general nonnegative values of $\rho$, it suffices to obtain \eqref{basicest}. For, in that case, the $d$ base points of any inclusion \eqref{seriesinseries} induce at least $d(s-d+r)-\rho$ shifts of vanishing order indices of the $g^{s-d+r}_m$. Just as before, the maximum possible number of shifts is $(s-d+r+1)(d-r)$, so we have
\[
d(s-d+r)- \rho \leq (s-d+r+1)(d-r),
\]
which contradicts $\rho+ \mu <0$.

\item The proof of Lemma 1 generalizes easily to a proof of \eqref{basicest} when $d^{\pr}=1$; as a consequence, we obtain a proof of the nonexistence of \eqref{seriesinseries} on the general curve when $\rho+\mu<0$ {\it under the additional assumption that the $d$ base points of the included $g^{s-d+r}_m$ are simple.}

\item The basic estimate \eqref{basicest} is strictly weaker than dimensional transversality for the Schubert varieties on $\mbox{Fl}(s-d+r,s; m)$ associated with the inclusion $g^{s-d+r}_{m-d^{\pr}}+ d^{\pr}p \hra g^s_m$ at 0, 1, $\infty$, and $p$, which predicts that
\beg{equation}\label{basicest2}
\sum_{k=0}^{s-d+r} k^{\pr}(j) \geq d^{\pr}(s-d+r+1)- \rho_Y(\al; \be; (0,1,\dots,1); s,m).
\end{equation}
Indeed, it is conceivable that dimensional transversality fails. In such an instance, the point $t=p$ belongs to the discriminant associated with the corresponding intersection of Schubert cycles $\Om_0$, $\Om_1$, $\Om_{\infty}$, and $\Om_t$. However, the discrepancy between \eqref{basicest} and \eqref{basicest2} suggests an alternative strategy for establishing \eqref{basicest}. Namely, given any inclusion $g^{s-d+r}_{m-d^{\pr}}+ d^{\pr}p \hra g^s_m$, fix a basis $(\wt{\sig}_j)_{j=0}^{s-d+r}$ for the $g^{s-d+r}_m$ whose orders of vanishing at 0 are strictly increasing with $j$, and whose orders of vanishing at $\infty$ are distinct. Omitting a single generator $\wt{\sig}_j$ from this basis determines a subseries $g^{s-d+r-1}_m$ of the $g^s_m$, again with a $d^{\pr}$-fold base point at $p$. As such, it belongs to the intersection of Schubert varieties $\wt{\Om}^{(j)}_0$, $\wt{\Om}^{(j)}_1$, $\wt{\Om}^{(j)}_{\infty}$, and $\wt{\Om}^{(j)}_t$, with $t=p$. The estimate \eqref{basicest} would follow provided we knew that for distinct choices $j_1, j_2$ of $j$, the sets of discriminantal values $t \in \mb{P}^1$-- those values of $t$ for which the intersections
\[
\wt{\Om}^{(j)}_0 \cap \wt{\Om}^{(j)}_1 \cap \wt{\Om}^{(j)}_{\infty} \cap \wt{\Om}^{(j)}_t, j=j_1, j_2
\]
fail to be transverse-- were mutually disjoint.
\end{itemize} 
We next prove a finiteness result for linear series with exceptional secant planes on a general curve in the case where $\rho=1$.

\beg{thm}\label{thm3}
If $\rho=1$ and $\mu=-1$, then there are finitely many linear series $g^s_m$ with $d$-secant $(d-r-1)$-planes on a general curve $C$ of genus $g$.
\end{thm}

\beg{proof}
Since the space of linear series on a general curve is irreducible whenever $\rho$ is positive, it suffices to show that some linear series without $d$-secant $(d-r-1)$-planes exists on $C$. To this end, it suffices to show that some smoothable linear series without $d$-secant $(d-r-1)$-planes exists on a flag curve $\widetilde{C}$ obtained by specialization from $C$.

We construct a particular choice of flag curve and linear series as follows. Fix a smooth irreducible elliptic curve $\widetilde{E}$ with general $j$-invariant, together with a general curve $\widetilde{Y}$ of genus $(g-1)$. Next, specialize $\widetilde{E}$ and $\widetilde{Y}$ to flag curves $E$ and $Y$. Glue $E$ and $Y$ transversely, letting $q$ denote their intersection. Let
\[
C^{\pr}:= Y \cup_q E.
\]
Furthermore, let $G^s_m(C^{\pr})$ denote the space of limit linear series along $C^{\pr}$, and let
\[
G^s_m(C^{\pr})_{(1,1,\dots,1,1)}
\]
denote the subspace of $G^s_m(C^{\pr})$ comprising limit linear series $V_Y$ for which 
\beg{equation}\label{ramif1}
\al(V_Y,q) \geq (1,1,\dots,1,1).
\end{equation}

The vanishing sequence corresponding to $(1,1,\dots,1,1)$ is $(1,2,3,\dots, s,s+1)$; by LS1, we deduce that 
\[
a(V_E,q) \geq (m-s-1,m-s,m-s+1, \dots, m-3,m-2,m-1),
\]
i.e., that
\beg{equation}\label{ramif2}
\al(V_E,q) \geq (m-s-1, \dots, m-s-1).
\end{equation}
Now let 
\[
r_Y= (1,\dots,1) \text{ and } r_E=(m-s-1,\dots,m-s-1).
\]
The modified Brill-Noether numbers $\rho(Y,(r_Y)_q)$ and $\rho(E,(r_E)_q)$, which compute the expected dimensions of the spaces of limit linear series along $Y$ and $E$ with ramification at $q$ prescribed by \eqref{ramif1} and \eqref{ramif2}, respectively, are
\[
\rho(Y,(r_Y)_q)= \rho(g-1,s,m)- (s+1) = \rho(g,s,m) +s -(s+1)=0
\]
and
\[
\rho(E,(r_E)_q)= \rho(1,s,m)- (s+1)(m-s-1)= 1.
\]

Since $\widetilde{Y}$ and $\widetilde{E}$ are general, their respective spaces of limit linear series $G^s_m(Y,(r_Y)_q)$ and $G^s_m(E,(r_E)_q)$ are of expected dimension, by Eisenbud and Harris' generalized Brill--Noether theorem \cite{EH3}. It follows immediately that $G^s_m(C^{\pr})_{(1,\dots,1)}$ is of expected dimension, so every linear series in $G^s_m(C^{\pr})_{(1,\dots,1)}$ smooths, by the Regeneration Theorem \cite[Thm 5.41]{HM}.

To prove Theorem~\ref{thm3}, it now suffices to show that no limit linear series in $G^s_m(C^{\pr})_{(1,\dots,1)}$ admits an inclusion \eqref{seriesinseries}. Note, however, that 
\[
a(V_Y,q) \geq (1,\dots,1)
\]
implies that along any component of the spine of $C^{\pr}$, any $g^s_m$ satisfies
\[
b_{s-i} \geq m-1-a_i
\]
for every index $i \in \{0,\dots,s\}$. (This is clear along $E$, where the special points $0$ and $\infty$ have vanishing sequences $(0,1,\dots,s)$ and $(m-s-1,m-s,\dots,m-1)$, and along $Y$ it follows from the fact that $\rho(Y,(r_Y)_q)=0$.) It now follows by the same argument used to prove Theorem 1 that no limit linear series in $G^s_m(C^{\pr})_{(1,\dots,1)}$ admits an inclusion \eqref{seriesinseries}.
\end{proof}


\beg{figure}
\beg{picture}(40,20)
\qbezier(0,3)(20,2)(58.5,15)
\qbezier(45,15)(67.5,4)(100,15)
\qbezier(70,15)(75,10)(70,-10) 
\qbezier(95,15)(107.5,4)(150,15)
\qbezier(120,15)(125,10)(120,-10) 
\qbezier(145,15)(167.5,4)(200,15)
\qbezier(170,15)(175,10)(170,-10) 
\put(185,17.5){$Y$}
\put(205,10){\dots}
\qbezier(222,15)(254.5,4)(270,15)
\qbezier(245,15)(250,10)(245,-10) 
\qbezier(265,15)(304.5,3)(327,10) 
\put(267,20){$q$}
\qbezier(305,15)(310,10)(305,-10) 
\put(302.5,17.5){$E$}
\put(125,-30){$\overbrace{(g-1) \text{ elliptic tails along }Y}$}
\end{picture}
\vspace{30pt}
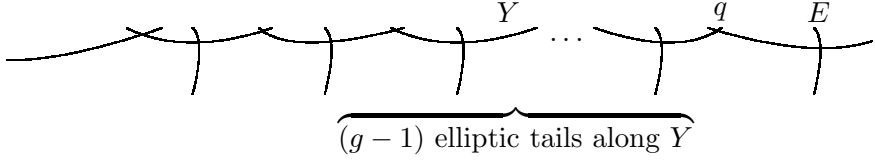
\caption{$C^{\pr}=Y \cup E$. Here $a(V_Y,q)= (1,\dots,1) \text{ and } a(V_E,q)=(m-s-1,\dots,m-s-1)$.}
\end{figure}

\section{Enumerative study}\label{enum}

In this section, we will study a special case of the following problem. Let $\pi: \mc{X} \ra B$ denote a one-parameter (flat) family of curves whose generic fiber is smooth, with some finite number of special fibers that are irreducible curves with nodes. We equip each fiber of $\pi$ with an $s$-dimensional series $g^s_m$. That is, $\mc{X}$ comes equipped with a line bundle $\mc{L}$, and on $B$ there is a vector bundle $\mc{V}$ of rank $(s+1)$, such that 
\[
\mc{V} \hra \pi_* \mc{L}.
\]
If $\mu=-1$, we expect finitely many fibers of $\pi$ to admit linear series with $d$-secant $(d-r-1)$-planes. We then ask for a formula for the number of such series, given in terms of tautological invariants associated with the family $\pi$.

One natural approach to the problem is to view those fibers whose associated linear series admit $d$-secant $(d-r-1)$-planes as a degeneracy locus for a map of vector bundles over $B$. This is the point of view adopted by Ziv Ran in his work \cite{R2,R3} on Hilbert schemes of families of nodal curves. Used in tandem with Porteous' formula for the class of a degeneracy locus of a map of vector bundles, Ran's work shows that the number of $d$-secant $(d-r-1)$-planes is a function $N^{d-r-1}_d$ of tautological invariants of the family $\pi$, namely:
\beg{equation}\label{invariants}
\al:=\pi_*(c_1^2(\mc{L})), \be:=\pi_*(c_1(\mc{L}) \cdot \om), \ga:=\pi_*(\om^2), \de_0, \text{ and } c:=c_1(\mc{V})
\end{equation}
where $\om= c_1(\om_{\mc{X}/B})$ and where $\de_0$ denotes the locus of points $b \in B$ for which the corresponding fiber $\mc{X}_b$ is singular.

In other words, for any fixed choice of $s$, we have
\vspace{-5pt}
\beg{equation}\label{basicformula}
N^{d-r-1}_d= P_{\al} \al+ P_{\be} \be+ P_{\ga} \ga+ P_c c + P_{\de_0} \de_0
\end{equation}
where the arguments $P$ are polynomials in $m$ and $g$ with coefficients in $\mb{Q}$. Unfortunately, the computational complexity of the calculus developed by Ran to evaluate $N^{d-r-1}_d$ grows exponentially with $d$. On the other hand, given that a formula \eqref{basicformula} in tautological invariants exists, the problem of evaluating it reduces to producing sufficiently many relations among the coefficients $P$.

In fact, the polynomials $P$ satisfy one ``obvious" relation. The most succinct explanation of how it arises involves ``twisting" our universal formula \eqref{basicformula} by a $\mb{Q}$-divisor. Namely, since $\mc{V} \hookrightarrow \pi_* \mc{L}$, the push-pull formula implies that for any $\mb{Q}$-divisor $D$ on $B$,
\beg{equation}\label{twistedinclusion}
\mc{V} \otimes \mc{O}(D) \hookrightarrow \pi_*(\mc{L} \otimes \pi^* \mc{O}(D)).
\end{equation}
Strictly speaking, the bundle to the right in \eqref{twistedinclusion} is only defined when $D$ is integral. However, we may obtain an integral version of \eqref{twistedinclusion} by passing to a finite cover of $B$. Doing so leaves \eqref{basicformula} unchanged. 

On the other hand, we have $c_1(\mc{V} \otimes \mc{O}(D))=0$ whenever $(s+1)D= \pi^* c_1(\mc{V})$. In that case, the renormalization
\[
\mc{L} \mapsto \mc{L} \otimes \pi^* \mc{O}(D)
\]
has the effect of trivializing $\mc{V}$.
The formula \eqref{basicformula} is clearly invariant under such renormalizations. Invariance implies that
\[
\beg{split}
&P_{\al} \pi_{*} \biggl(c_1(\mc{L})- \fr{\pi^{*}c}{s+1}\biggr)^2+ P_{\be} \pi_*\biggl( \biggl(c_1(\mc{L})- \fr{\pi^*c}{s+1}\biggr) \cdot \om \biggr) + P_{\ga} \ga+ P_{\de_0} \de_0 \\
&=P_{\al} \pi_*(c_1^2(\mc{L}))+ P_{\be} \pi_*(c_1(\mc{L}) \cdot \om) + P_{\ga} \ga+ P_{\de_0} \de_0+ P_c c.
\end{split}
\]
where as usual $c= c_1(\mc{V})$. The coefficient of $c$ in the left-hand expression is $-\fr{2m}{s+1} P_{\al}- \fr{2g-2}{s+1} P_{\be}$; since the coefficient of $c$ on the right-hand expression is $P_c$, we deduce that
\beg{equation}\label{obveqn}
2m P_{\al}+ (2g-2) P_{\be}+ (s+1) P_c=0.
\end{equation}

Note that Theorem 1 implies that when $\rho=0$ and $\mu=-1$, curves that admit linear series with exceptional secant planes sweep out a divisor $\mbox{Sec}$ in $\ov{\mc{M}}_g$. Given the polynomials $P$, the class of $\mbox{Sec}$, modulo higher-boundary divisors $\de_i, i \geq 1$, may be determined via Khosla's pushforward formulas in \cite{Kh}. All of this is explained in detail in \cite{Co1} and \cite{Co2}, where the slopes of secant plane divisors are computed.

Hereafter in this paper, where our focus is linear series on the general curve, we will assume that $\pi$ is a trivial family; more specifically, that $\mc{X}=B \times C$ where $C$ is a smooth curve. We then have $\ga= \de_0=0$, so only two more relations are required to determine the tautological coefficients $P$. For this purpose, we evaluate our secant-plane formula \eqref{basicformula} along test families involving projections of a fixed curve in projective space.

\subsection{Test families}\label{testfamilies}
Our test families are as follows:
\vspace{-4pt}
\beg{enumerate}
\item {\it Family one.} Projections of a generic curve of degree $m$ in $\mb{P}^{s+1}$ from points along a disjoint line. 
\item {\it Family two.} Projections of a generic curve of degree $m+1$ in $\mb{P}^{s+1}$ from points along the curve.
\end{enumerate}
\vspace{-4pt}
Now assume that $\mu(d,r,s)=-1$. Let $A$ denote the expected number of $d$-secant $(d-r)$-planes to a curve of degree $m$ and genus $g$ in $\mb{P}^{s+1}$ that intersect a general line. Let $A^{\pr}$ denote the expected number of $(d+1)$-secant $(d-r)$-planes to a curve of degree $(m+1)$ and genus $g$ in $\mb{P}^{s+1}$.
The expected number of fibers of the first (resp., second) family with $d$-secant $(d-r-1)$-planes equals $A$ (resp., $(d+1)A^{\pr}$). 

Determining those relations among the tautological coefficients induced by the two families requires knowing the values of $\al,  \be, \text{ and }\ga$ along each family $\pi: \mc{X} \ra B$.
These are determined as follows.
\beg{itemize}
\item {\it Family one.} The base and total spaces of our family are $B=\mb{P}^1$ and $\mc{X}= \mb{P}^1 \times C$, respectively. Letting $\pi_1$ and $\pi_2$ denote, respectively, the projections of $\mc{X}$ onto $\mb{P}^1$ and $C$, we have
\[
\mc{L}= \pi_2^* \mc{O}_C(1), \om_{\mc{X}/\mb{P}^1}= \pi_2^* \om_C, \text{ and } \mc{V}= \mc{O}_{\mb{G}}(-1) \otimes \mc{O}_{\mb{P}^1}
\]
where $\mb{G}=\mb{G}(s,s+1)$ denotes the Grassmannian of hyperplanes in $\mb{P}^{s+1}$. Accordingly,
\[
\al= \be= \ga= \de_0= 0, \text{ and } c=-1.
\]
It follows that
\[
P_c= -A.
\]

\item {\it Family two.} This time, $\mc{X}= C \times C$ and $B=C$. Here
\[
\mc{L}= \pi_2^* \mc{O}_C(1) \otimes \mc{O}(-\Delta), \om_{\mc{X}/\mb{P}^1}= \pi_2^* \om_C, \text{ and } \mc{V}= \mc{O}_{\mb{G}}(-1) \otimes \mc{O}_C.
\]
Consequently, letting $H=c_1(\mc{O}_C(1))$, we have
\[
\beg{split}
\al &=-2 \Delta \cdot \pi_2^* (m+1)\{\mbox{pt}_C\}+
\Delta^2 = -2m-2g, \\
\be &= (\pi_2^*H -\Delta) \cdot \pi_2^*K_C = 2-2g, \\
c &= -m-1, \text{ and } \ga = \de_0= 0.
\end{split}
\]
It follows that
\[
(-2m-2g)P_{\al}+ (2-2g)P_{\be}+ (-m-1)P_c= (d+1)A^{\pr}.
\]
\end{itemize}

\subsection{Classical formulas for $A$ and $A^{\pr}$, and their significance}\label{AFormula}

Formulas for $A$ and $A^{\pr}$ were calculated by Macdonald \cite{M} and Arbarello, et. al, in \cite[Ch. 8]{ACGH}. The formulas have enumerative significance only when the loci in question are actually zero-dimensional. On the other hand, for the purpose of calculating class formulas for secant-plane divisors on $\ov{\mc{M}}_g$, it clearly suffices to show that for every fixed triple $(d,r,s)$, Macdonald's formulas are enumerative whenever $m=m(d,r,s)$ is sufficiently large.

To do so, we view the curve $C \sub \mb{P}^{s+1}$ in question as the image under projection of a {\it non-special} curve $\widetilde{C}$ in a higher-dimensional ambient space. We then re-interpret the secant behavior of $C$ in terms of the secant behavior of $\widetilde{C}$; the latter, in turn, may be characterized completely because $\widetilde{C}$ is non-special.

Given a curve $\widetilde{C}$, let $L$ be a line bundle of degree $\widetilde{m}$ on $\widetilde{C}$, let $V \sub H^0(\widetilde{C},L)$; the pair $(L,V)$ defines a linear series on $\widetilde{C}$. Now let $T^{\wt{d}}(L)$ denote the vector bundle
\[
T^{\wt{d}}(L)= (\pi_{1,\dots,\widetilde{d}})_* (\pi_{\widetilde{d}+1}^* L \otimes \mc{O}_{\wt{C}^{\widetilde{d}+1}}/\mc{O}_{\wt{C}^{\widetilde{d}+1}}(-\De_{\widetilde{d}+1}))
\]
over $\wt{C}^{\widetilde{d}}$, where $\pi_i, i=1 \dots \widetilde{d}+1$ denote the $\widetilde{d}+1$ projections of $\wt{C}^{\widetilde{d}+1}$ to $\widetilde{C}$, $\pi_{1,\dots,\widetilde{d}}$ denotes the product of the first $\widetilde{d}$ projections, and $\De_{\widetilde{d}+1} \sub \wt{C}^{\widetilde{d}+1} $ denotes the ``big" diagonal of $(\widetilde{d}+1)$-tuples whose $i$th and $(\widetilde{d}+1)$st coordinates are the same. The bundle $T^{\wt{d}}(L)$ has fiber $H^0(L/L(-D))$ over a divisor $D \sub \wt{C}^{\widetilde{d}}$.

Note that  the $\widetilde{d}$-secant $(\widetilde{d}-\widetilde{r}-1)$-planes to the image of $\widetilde{C}$ under $(L,V)$ correspond to the sublocus of $\wt{C}^{\widetilde{d}}$ over which the evaluation map
\beg{equation}\label{ev}
V \st{\mbox{ev}}{\lra} T^{\wt{d}}(L)
\end{equation}
has rank $(\widetilde{d}-\widetilde{r})$. 

Moreover, by Serre duality,
\beg{equation}\label{serredual}
H^0(\om_{\widetilde{C}} \otimes L^{\vee} \otimes \mc{O}_{\widetilde{C}}(p_1 + \dots + p_{\widetilde{d}}))^{\vee} \cong H^1(L(-p_1- \dots -p_{\widetilde{d}}));
\end{equation}
both vector spaces are zero whenever $\om_{\widetilde{C}} \otimes L^{\vee} \otimes \mc{O}_{\widetilde{C}}(p_1 + \dots + p_{\widetilde{d}})$ has negative degree. In particular, whenever
\beg{equation}\label{deghypoth}
\widetilde{m} \geq 2g-1+ \widetilde{d},
\end{equation}
the vector space on the right-hand side of \eqref{serredual} is zero. It follows that the evaluation map \eqref{ev} is surjective for the complete linear series $(L,H^0(\mc{O}_{\widetilde{C}}(D))$ whenever $D \sub \widetilde{C}$ is a divisor of degree $\widetilde{m}$ verifying \eqref{deghypoth}. Equivalently, whenever \eqref{deghypoth} holds, every $\widetilde{d}$-tuple of points in $\widetilde{C}$ determines a secant plane to the image of $(L,H^0(\mc{O}_{\widetilde{C}}(D))$ is of maximal dimension $(\widetilde{d}-1)$.

Now let $\widetilde{s}:= h^0(\mc{O}_{\widetilde{C}}(D))$. Somewhat abusively, we will identify $\widetilde{C}$ with its image in $\mb{P}^{\widetilde{s}}$. Let $C$ denote the image of $\widetilde{C}$ under projection from an $(\widetilde{s}-s-2)$-dimensional center $\Ga \sub \mb{P}^{\widetilde{s}}$ disjoint from $\widetilde{C}$.

Note that $\widetilde{d}$-secant $(\widetilde{d}-\widetilde{r}-1)$-planes to $C$ are in bijective correspondence with those $\widetilde{d}$-secant $(\widetilde{d}-1)$-planes to $\widetilde{C}$ that have at least $(\widetilde{r}-1)$-dimensional intersections with $\Ga$. These, in turn, comprise a subset $\mc{S} \sub \mb{G}(\widetilde{d}-1,\widetilde{s})$ defined by
\beg{equation}\label{intersection}
\mc{S}= \mc{V} \cap \sig_{\ubr{s-\widetilde{d}+\widetilde{r}+2,\dots,s-\widetilde{d}+\widetilde{r}+2}_{\widetilde{r} \text{ times}}}
\end{equation}
where $\mc{V}$, the image of $\wt{C}^{\widetilde{d}}$ in $\mb{G}(\widetilde{d}-1,\widetilde{s})$, is the variety of $\widetilde{d}$-secant $(\widetilde{d}-1)$-planes to $\widetilde{C}$, and
the term involving $\sig$ denotes the Schubert cycle of $(\widetilde{d}-1)$-planes to $\widetilde{C}$ that have at least $(\widetilde{r}-1)$-dimensional intersections with $\Ga$. For a general choice of projection center $\Ga$, the intersection \eqref{intersection} is transverse; it follows that
\beg{equation}\label{expdim}
\dim \mc{S}= \widetilde{d}- \widetilde{r}(s- \widetilde{d}+ \widetilde{r}+2),
\end{equation}
In particular, if $\widetilde{d}=d+1$ and $\widetilde{r}=r$, then $\dim \mc{S}= 1+ \mu(d,r,s)= 0$, which shows that for any choice of $(d,r,s)$, the formula for $A^{\pr}$ is enumerative whenever $m=m(d,r,s)$ is chosen to be sufficiently large.

Similarly, to handle $A$, note that there is a bijection between $\widetilde{d}$-secant $(\widetilde{d}-\widetilde{r}-1)$-planes to $C$ that intersect a general line and $\widetilde{d}$-secant $(\widetilde{d}-1)$-planes to $\widetilde{C}$ that have at least $(\widetilde{r}-1)$-dimensional intersections with $\Ga$, and which further intersect a general line $l \sub \mb{P}^{\widetilde{s}}$. These, in turn, comprise a subset $\mc{S}^{\pr} \sub \mb{G}(\widetilde{d}-1,\widetilde{s})$ given by
\beg{equation}\label{intersection2}
\mc{S}^{\pr}= \mc{V} \cap \sig_{\ubr{s-\widetilde{d}+\widetilde{r}+2,\dots,s-\widetilde{d}+\widetilde{r}+2}_{\widetilde{r} \text{ times}},s-\widetilde{d}+\widetilde{r}+1}.
\end{equation}
For a general choice of projection center $\Ga$ and line $l$, the intersection \eqref{intersection2} is transverse. In particular, if $\widetilde{d}=d$ and $\widetilde{r}=r-1$, then $\dim \mc{S}^{\pr}= 0$, which shows that for any choice of $(d,r,s)$, the formula for $A$ is enumerative whenever $m$ is sufficiently large.

\subsection{The case $r=1$}\label{Ndfmla} 
Note that the equation $\mu=-\al-1$ may be rewritten in the following form:
\[
s=\fr{d+\al+1}{r}+ d-1-r.
\]
As a result, $r$ necessarily divides $(d+\al+1)$, say $d=\ga r-\al-1$, and correspondingly,
\[
s=(\ga-1)r+ \ga-\al-2.
\]
In particular, whenever $\rho=0$ and $\mu=-1$, we have $1 \leq r \leq s$. As a result, we will focus mainly on the two ``extremal" cases of series where $r=1$ or $r=s$. Our strongest results are for $r=1$; accordingly, we treat this case in this subsection and the two following it.

As a special case of \cite[Ch. VIII, Prop. 4.2]{ACGH}, the expected number of $(d+1)$-secant $(d-1)$-planes to a curve $C$ of degree $(m+1)$ and genus $g$ in $\mb{P}^{2d}$ is
\begin{equation}\label{ACGHNd}
A^{\pr}= \sum_{\al=0}^{d+1} (-1)^{\al} \binom{g+2d-(m+1)}{\al} \binom{g}{d+1-\al}.
\end{equation}
In fact, the formula for $A$ in case $r=1$ is implied by the preceding formula. To see why, note that $d$-secant $(d-1)$-planes to a curve $C$ of degree $m$ and genus $g$ in $\mb{P}^{2d}$ that intersect a disjoint line $l$ are in bijection with $d$-secant $(d-2)$-planes to a curve $C$ of degree $m$ and genus $g$ in $\mb{P}^{2d-2}$ (simply project with center $l$). It follows that
\[
A= \sum_{\al=0}^d (-1)^{\al} \binom{g+2d-(m+3)}{\al} \binom{g}{d-\al}.
\]
{\bf Remark.} Denote the generating function for the formulas $A=A(d,g,m)$ in case $r=1$ by $\sum_{d \geq 0} N_d(g,m) z^d$, where
\[
N_d(g,m):= \# \text{ of }d-\text{secant }(d-2)-\text{planes to a }g^{2d-2}_m \text{ on a genus-}g \text{ curve}. 
\]
(As a matter of convention, we let $N_0(g,m)=1$, and $N_1(g,m)=d$.)

The generating function for $N_d(g,m)$ is as follows (here we view $g$ and $m$ as fixed, and we allow the parameter $d$ to vary).
\beg{thm}\label{Ndthm}
\beg{equation}\label{N_dbis}
\sum_{d \geq 0} N_d(g,m) z^d= \biggl(\fr{2}{(1+4z)^{1/2}+1} \biggr)^{2g-2-m} \cdot (1+4z)^{\fr{g-1}{2}}.
\end{equation}
\end{thm}

\beg{proof}
We will in fact prove that
{\small
\beg{equation}\label{N_d}
\sum_{d \geq 0} N_d(g,m) z^d= \exp \biggl({\sum_{n>0} \fr{(-1)^{n-1}}{n} \biggl[\binom{2n-1}{n-1}m+ \biggl(4^{n-1}- \binom{2n-1}{n-1}\biggr) (2g-2)\biggr] z^n}\biggr). 
\end{equation}
}
To see that the formulas \eqref{N_d} and \eqref{N_dbis} are equivalent, begin by recalling that the generating function $C(z)= \sum_{n \geq 0} C_n z^n$ for the Catalan numbers $C_n= \fr{\binom{2n}{n}}{n+1}$ is given explicitly by
\[
C(z)= \fr{1- \sqrt{1-4z}}{2z}.
\]
On the other hand, we have
\[
\fr{\binom{2n-1}{n-1}}{n}= \biggl(2- \fr{1}{n}\biggr) C_{n-1}.
\]
See, e.g., \cite[Sect 2]{De1} for generalities concerning Catalan numbers. It follows that \eqref{N_d} may be rewritten as follows:
{\small
\[
\beg{split}
\sum_{d \geq 0} N_d(g,m) z^d&= \exp \biggl[{\sum_{n>0} (-1)^{n-1} \biggl[ \biggl[\biggl(2-\fr{1}{n}\biggr)(m-2g+2) C_{n-1} z^n} \biggr]+ 4^{n-1} \cdot (2g-2) \fr{z^n}{n} \biggr] \biggr]\\
&= \exp \biggl[(2m-4g+4)z C(-z)- (m-2g+2) \int C(-z) dz + (2g-2) \int \fr{1}{1+4z} dz \biggr].
\end{split}
\]
}
Here $\int$ denotes integration of formal power series according to the convention that
\[
\int z^n dz= \fr{1}{n+1} z^{n+1}
\]
for all nonnegative integers $n$. 
We now claim that
{\small
\beg{equation}\label{id1}
\beg{split}
-\int C(-z) dz&= \int \fr{1-(1+4z)^{1/2}}{2z} dz \\
&= 1-(1+4z)^{1/2}+\fr{1}{2} \ln \fr{z}{(1+4z)^{1/2}-1}+ \fr{1}{2}\ln((1+4z)^{1/2}+1).
\end{split}
\end{equation}
}
Indeed, it is easy to check that the sum of the derivatives of the power series on the right side equals $-\fr{2}{(1+4z)^{1/2}+1}$, or equivalently, $\fr{1-(1+4z)^{1/2}}{2z}$. Moreover, l'H\^{o}pital's rule shows that $\ln \fr{z}{(1+4z)^{1/2}-1}$ evaluates to $-\ln(2)$ at $z=0$, so the right side of \eqref{id1} evaluates to zero at $z=0$, as required according to our convention regarding integration. The claim follows.

As a result, we may write
{\small
\[
\beg{split}
\sum_{d \geq 0} N_d(g,m) z^d 
&= \exp \biggl[(2g-2-m) \biggl( \fr{1}{2} \ln \fr{z}{(1+4z)^{1/2}-1}+ \fr{1}{2}\ln((1+4z)^{1/2}+1)\bigg)+ \fr{(g-1)}{2} \ln(1+4z) \biggr].
\end{split}
\]
}
Since
\[
\fr{\ln((1+4z)^{1/2}-1))}{2}+ \fr{\ln((1+4z)^{1/2}+1))}{2}= \fr{\ln 4z}{2}= \fr{\ln z}{2}+ \ln 2,
\]
we have
\[
\sum_{d \geq 0} N_d(g,m) z^d = \exp \biggl[(2g-2-m)\biggl(\ln 2- \ln((1+4z)^{1/2}+1)\biggr) + \fr{(g-1)\ln (1+4z)}{2} \biggr],
\]
from which \eqref{N_dbis} follows immediately.

\medskip
To prove \eqref{N_d}, proceed as follows. Begin by fixing a positive integer $d>0$, and let $C$ denote the image of a $g^{2d-2}_m$ that is sufficiently ``nonspecial" in the sense of the preceding section. Then, as noted in the preceding section, $N_d(g,m)$ computes the degree of the locus of $d$-tuples in $\mbox{Sym}^d C$ for which the evaluation map \eqref{ev} has rank $(d-1)$. In fact, we will find it more convenient to work on the $d$-tuple product $C^d$. Clearly, $N_d(g,m)$ computes $\fr{1}{d!}$ times the degree $\widetilde{N_d}(g,m)$ of the locus along which the corresponding evaluation map has rank $(d-1)$, since there are $d!$ permutations of any given $d$-tuple corresponding to a given $d$-secant plane.

On the other hand, Porteous' formula implies that $\widetilde{N_d}(g,m)$ is equal to the degree of the determinant
\beg{equation}\label{det}
\left| \beg{array}{ccccc}
c_1 & c_2 & \cdots  & c_{d-1} & c_d \\
1 & c_1 & \cdots  & c_{d-2} & c_{d-1} \\
\cdots & \cdots & \cdots & \cdots & \cdots \\
0 & \cdots & 0 & 1 & c_1
\end{array} \right|
\end{equation}
where $c_i$ denotes the $i$th Chern class of the secant bundle $T^d(L)$ over $C^d$. Note \cite{R1} that the Chern polynomial of $T^d(L)$ is given by
\[
c_t(T^d(L))= (1+ l_1 t) \cdot (1+ (l_2-\De_2)t) \cdots (1+ (l_d-\De_d)t)
\]
where $l_i, 1 \leq i \leq d$ is the pullback of $c_1(L)$ along the $i$th projection $C^d \ra C$, and $\De_j, 2 \leq j \leq d$ is the (class of the) diagonal defined by
\[
\De_j= \{(x_1,\dots,x_d) \in C^d | x_i=x_j \text{ for some } i<j \}.
\]
In particular, modulo $l_i$'s, we have
\[
c_i= (-1)^i s_i(\De_2,\dots, \De_d)
\]
where $s_i$ denotes the $i$th elementary symmetric function. 

The degree of \eqref{det} may be expressed as a polynomial in $m$ and $(2g-2)$ with integer coefficients. To see why, note that
\beg{equation}\label{diag}
\De_j= \sum_{i=1}^{j-1} \De_{i,j}
\end{equation}
for every $2 \leq j \leq d$. Furthermore, letting $p_i$ denote the projection of $C^d$ to the $i$th copy of $C$, we have
\[
l_j \cdot \De_{i,j}= p_i^* m \{ \mbox{pt}_C \},
\]
and
\[
\De_{i,j}^2= -p_i^* \om_C \cdot \De_{i,j}= -(2g-2)p_i^* \{\mbox{pt}_C\} \cdot \De_{i,j}
\]
for every choice of $(i,j)$. Consequently, \eqref{det} is a sum of degree-$d$ monomials in the diagonal summands $\De_{i,j}$ and the $l_k$. 

Now, letting $s_i(x_1,\dots,x_d)$ denote the $i$th elementary symmetric function in the indeterminates $x_i$, an easy inductive argument shows that for every positive integer $n$,
\beg{equation}\label{matrixid}
\left| \beg{array}{ccccc}
s_1 & s_2 & \cdots & s_{n-1} & s_n \\
1 & s_1 & \cdots & s_{n-2} & s_{n-1} \\
\cdots & \cdots & \cdots & \cdots & \cdots \\
0 & \cdots & 0 & 1 & s_1
\end{array} \right| = \mathop{\sum_{i_1, \dots, i_{d} \geq 0}}_{i_1+ \cdots i_{d}=n} x_1^{i_1} \cdots x_d^{i_{d}}.
\end{equation}


The identity \eqref{matrixid} implies that the term of degree one in $(2g-2)$ and zero in $m$ of the determinant \eqref{det} is equal to the term of appropriate degree in
\beg{equation}\label{const}
(-1)^d \mathop{\sum_{i_1, \dots, i_{d-1} \geq 0}}_{i_1+ \cdots i_{d-1}=d} \De_2^{i_1} \cdots \De_d^{i_{d-1}}.
\end{equation}

Similarly, the term of degree zero in $(2g-2)$ and one in $m$ of \eqref{det} is equal to the term of corresponding degree in
\beg{equation}\label{mterm}
(-1)^{d-1} \mathop{\sum_{i_1, \dots, i_{d-1} \geq 0}}_{i_1+ \cdots i_{d-1}=d-1} \sum_{j=1}^{d} a_j l_j \De_2^{i_1} \cdots \De_d^{i_{d-1}}
\end{equation}
where $a_j=1$ if $j=1$ and $a_j=i_j+1$ whenever $2 \leq j \leq d$.

On the other hand, as an immediate consequence of the way in which the coefficients $a_j$ are defined, the intersection \eqref{mterm} pushes down to
\beg{equation}\label{mterm2}
\beg{split}
&(-1)^{d-1} \mathop{\sum_{i_1, \dots, i_{d-1} \geq 0}}_{i_1+ \cdots i_{d-1}=d-1} \biggl(1+ \sum_{j=1}^{d-1}(i_j+1)\biggr) \De_2^{i_1} \cdots \De_d^{i_{d-1}} \\
&= (-1)^{d-1} (2d-1) \mathop{\sum_{i_1, \dots, i_{d-1} \geq 0}}_{i_1+ \cdots i_{d-1}=d-1} \De_2^{i_1} \cdots \De_d^{i_{d-1}}.
\end{split}
\end{equation}

\beg{lem}
Up to a sign, the term of degree zero in $(2g-2)$ and degree one in $m$ in \eqref{mterm2} is equal to 
\[
\binom{2d-1}{d-1} (d-1)! \cdot m.
\]
\end{lem}

\beg{lem}
Up to a sign, the term of degree one in $(2g-2)$ and zero in $m$ in \eqref{const} is equal to 
\[
\biggl(4^{d-1}- \binom{2d-1}{d-1}\biggr) (d-1)! \cdot (2g-2).
\]
\end{lem}

To go further, the following observation will play a crucial r\^{o}le. For any $d \geq 2$, let $K_d$ denote the complete graph on $d$ labeled vertices $v_1, \dots, v_d$, whose edges $e_{i,j}=\ov{v_i v_j}$ are each oriented with arrows pointing towards $v_j$ whenever $i<j$. The degree of our determinant \eqref{det} computes a sum of monomials involving $\De_i$ and $l_j$, where $2 \leq i \leq d$ and $1 \leq j \leq d$, and so may be viewed as a tally $\widetilde{S}$ of (not-necessarily connected) subgraphs of $K_d$, each counted with the appropriate weights. By the Exponential Formula \cite[5.1.6]{St}, the exponential generating function for the latter, as $d$ varies, is equal to $e^{E_S}$, where $E_S$ is the exponential generating function for the corresponding tally of connected subgraphs, which correspond, in turn, to the intersections described in Lemmas 4 and 5.

More precisely now, fix an integer $d \geq 2$, and consider subgraphs of $K_d$ having some number $\tau$ of connected components $G_1, \dots, G_{\tau}$. (Strictly speaking, we are not merely interested in subgraphs, but in graphs supported on $K_d$ in which at most one edge appears with multiplicity 2, so our terminology is abusive.) Say that the component subgraph $G_i$ has $n_e(i)$ vertices; we stipulate that either these are connected by $n_e(i)$ edges, or else that $G_i$ has a unique ``marked" vertex and $(n_e(i)-1)$ edges.   Marked vertices $v_j$ correspond to instances of $l_j$, while edges $e_{i,j}$ correspond to small diagonals $\De_{i,j}= \{(x_1,\dots,x_d) \in C^d | x_i=x_j\}$ associated to $d$-tuples whose $i$th and $j$th coordinates agree. In the case where no marked vertex appears, at most one edge $e_{i,j}$ may appear with multiplicity 2, in which case it corresponds to $\De_{i,j}^2$. 

When $G_i$ has no marked vertices, assign to each vertex $v_j$ in $G_i$ a weight
\[
w_{G_i,j}= \binom{\mbox{indeg}(G_i,j)}{i_1,\dots,i_{j-1}}
\]
where $\mbox{indeg}(G_i,j)$ is equal to the indegree of $v_j$ in $G_i$, i.e., the total number of edges of $G_i$ incident with $v_j$, counted with their nonnegative multiplicities $i_1, \dots, i_{j-1}$. Let 
\[
w_{G_i}= \prod_j w_{G_i,j}
\]
where the product is over all vertices $v_j$ appearing in $G_i$.

Similarly, when $G_i$ contains a marked vertex, assign to each vertex $v_j$ in $G_i$ (including the marked vertex) the weight
\beg{equation}\label{wt1}
w_{G_i,j}= (\mbox{indeg}(G_i,j))!,
\end{equation}
and let
\beg{equation}\label{wt2}
w_{G_i}= (2n(e_i)+1) \prod_j w_{G_i,j}
\end{equation}
where the product is over all vertices $v_j$ appearing in $G_i$.

Now let
\[
P_{G_i}^{(1)}:= (-1)^{n(e_i)+1} w(G_i) (2g-2), \text{ and }
P_{G_i}^{(2)}:= (-1)^{n(e_i)} w(G_i) m.
\]

Set $P_G^{(k)}:= \prod_{i=1}^{\tau} P_{G_i}^{(k)}, k=1,2$. The $P_{G_i}^{(k)}, k=1,2$ correspond to monomial intersection products of the forms
\[
\De_{i_1,i_1^{\pr}} \cdots \De_{i_{n(e_i)},i_{n(e_i)}^{\pr}}, \text{ and } l_j \De_{i_1,i_1^{\pr}} \cdots \De_{i_{n(e_i)-1},i_{n(e_i)-1}^{\pr}},
\]
respectively, counted with weights $w(G_i)$ prescribed by \eqref{wt1} and \eqref{wt2}.
{\it A crucial point is that $P_G:=P_G^{(1)}+ P_G^{(2)}$ represents the contribution of the intersection product corresponding to $G$ to the degree of the determinant \eqref{det}}. In the case where $G$ is connected, this follows from systematically applying the diagonal class decomposition \eqref{diag} and expanding \eqref{const} and \eqref{mterm2} accordingly.  The collection of indegrees $\mbox{indeg}(G,j), 2 \leq j \leq d$ specifies a monomial either of type $\De_2^{i_1} \cdots \De_d^{i_{d-1}}$ or of type $l_j \De_2^{i_1} \cdots \De_d^{i_{d-1}}$ (depending upon whether the sum of indegrees is $d$ or $(d-1)$), and the weight $w_{G_i,j}$ is a binomial coefficient attached to $\De_{i,j}^{i_{j-1}}$ that appears when we expand \eqref{const} or \eqref{mterm2}, as the case may be. The generalization to the case where $G$ has multiple connected components is immediate.

On the other hand, it is not hard to see that given any subset $B \sub \{1,\dots d\}$, the values of the functions $f_1$ and $f_2$ that compute the weighted tallies of all connected subgraphs of the complete graph on $B$ with or without marked vertices, respectively, depend only on the cardinality of $B$. Let $\widetilde{f}_1$ and $\widetilde{f}_2$ denote the functions that compute the corresponding ``disconnected" weighted tallies of subgraphs of $K_d$. Allowing $d$ to vary, we obtain exponential generating functions $E_{f_i}$ and $E_{\widetilde{f}_i}$ for $f_i$ and $\widetilde{f}_i$, respectively, where $i=1, 2$. The Exponential Formula implies that $E_{f_i}$ and $E_{\widetilde{f}_i}$ are related by
\[
E_{\widetilde{f}_i}= \mbox{exp}(E_{f_i}),
\]
for $i=1, 2$.
Now let
\[
\widetilde{f}:= \sum_G P_G.
\]
Since every subgraph of $K_d$ of interest to us can be realized as the union of a subgraph (possibly disconnected) with marked vertices and a subgraph without marked vertices, the exponential generating function $E_{\widetilde{f}}$ of $\widetilde{f}$ satisfies
\[
E_{\widetilde{f}}= E_{\widetilde{f}_1} \cdot E_{\widetilde{f}_2}
\]
by \cite[Prop. 5.1.3]{St}.

\medskip
Consequently, to prove Theorem 4, we are reduced to proving Lemmas 4 and 5.

\medskip
{\fl \bf Proof of Lemma 4:}
The only terms of degree one in $m$ and zero in $(2g-2)$ in \eqref{mterm2} correspond to $(d-1)$-tuples $(i_1,\dots,i_{d-1})$ that satisfy the additional constraint
\beg{equation}\label{constraint}
\sum_{k=1}^j i_k \leq j, \text{ for all } 1 \leq j \leq d-1.
\end{equation}

Notice that the number of such $(d-1)$-tuples is exactly the $(d-1)$st Catalan number $C(d-1)$.

We now expand \eqref{mterm2} according to \eqref{diag}. The monomials of relevance in the resulting expanded intersection product are exactly those in which no diagonal factor $\De_{i,j}$ is repeated. 

Accordingly, proving the lemma now transposes into the following combinatorial problem. Let $K_d$ denote the complete graph on $d$ labeled vertices $v_1, \dots, v_d$, whose edges $e_{i,j}$ are marked as before. Consider the set $\mc{T}$ of connected spanning trees on $K_d$. To each vertex $v_j, 2 \leq j \leq d$ of a graph $G$ in $\mc{T}$, assign the weight
\[
w_{G,j}= (\mbox{indeg}(G,j))!.
\]
where $\mbox{indeg}(G,j)$ denotes the total indegree of the vertex $v_j$ in $G$. Now set $w_G= \prod_{2 \leq j \leq d} w_{G,j}$.
In light of \eqref{mterm2}, it then suffices to show that
\[
(2d-1) \sum_{G \in \mc{T}} w_G= \binom{2d-1}{d-1}(d-1)!,
\]
i.e., that
\beg{equation}\label{reduction}
\sum_{G \in \mc{T}} w_G= \fr{(2d-2)!}{d!}.
\end{equation}
Since $\mc{T}$ has $C(d-1)$ elements, \eqref{reduction} will follow provided we can show that the average weight $w_G$ over all $G$ in $\mc{T}$ equals $(d-1)!$.

To this end, let $a_{i_1,\dots,i_{d-1}}$ denote the number of connected spanning trees on $K_d$ with indegrees $i_1,\dots,i_{d-1}$ at vertices $v_2,\dots,v_d$. Clearly, we have
\[
\sum_{G \in \mc{T}} w_G= \sum_{i_1,\dots,i_{d-1}} a_{i_1,\dots,i_{d-1}} i_1! \dots i_{d-1}!
\]
where the $i_j, 1 \leq j \leq d-1$ are nonnegative integers whose sum equals $(d-1)$, and which satisfy the constraint \eqref{constraint}. It then suffices to show that for any given choice of $(d-1)$-tuple $(i_1,\dots,i_{d-1})$ satisfying our constraints, the average value of all $a_{j_1,\dots,j_{d-1}}$ arising from permuting $(i_1,\dots,i_{d-1})$ (while still respecting \eqref{constraint}) equals $\fr{(d-1)!}{ i_1! \dots i_{d-1}!}$. 

As a matter of terminology, let an {\it admissible} permutation of a given $(d-1)$-tuple $(i_1,\dots,i_{d-1})$ denote a $(d-1)$-tuple obtained by permuting $(i_1,\dots,i_{d-1})$ that satisfies \eqref{constraint}. Let $\phi(i_1,\dots,i_{d-1})$ denote the number of admissible permutations of a given $(d-1)$-tuple $(i_1,\dots,i_{d-1})$. Note that $\phi(i_1,\dots,i_{d-1})$ is exactly the number of Dyck paths of semilength $(d-1)$ associated to the corresponding partition $(\la_1^{e_1},\dots,\la_1^{e_l})$ of $(d-1)$, obtained by discarding every instance of zero in $(i_1,\dots,i_{d-1})$. (A {\it Dyck path of semilength $d$} is a sequence of successive symbols $U$ and $D$ (``up" and ``down") of length $2d$, with the property that at any given position $1 \leq j \leq 2d$, no more $D$'s than $U$'s lie to the left. See \cite[Sect 2]{De1}.) Here $\la_i^{e_i}$ denotes a sequence of $e_i$ identical terms $\la_i$. We then have
\beg{equation}\label{phieq}
\phi(i_1,\dots,i_{d-1})= \fr{(d-1)!}{(d-k)!e_1! \dots e_l!}
\end{equation}
where $k=\sum_{i=1}^l e_i$, by \cite[Thm. 5.3.10]{St}. Accordingly, it suffices to show that
\beg{equation}\label{Syma}
a_{\la}=\fr{(d-1)!}{(d-k)!e_1! \dots e_l!} \cdot \fr{(d-1)!}{ (\la_1!)^{e_1} \dots (\la_l!)^{e_l}}
\end{equation}
where $a_{\la}:= \sum_{(j_1,\dots,j_{d-1})} a_{(j_1,\dots,j_{d-1})}$ is the total number of connected spanning trees with indegree sequences $(j_1,\dots,j_{d-1})$ that are admissible permutations of a fixed indegree sequence $(i_1,\dots,i_{d-1})$ corresponding to the partition $\la=(\la_1^{e_1},\dots,\la_l^{e_l})$. But this is the main result of \cite{DY}. \qed

\medskip
Thus far we have been unable to obtain a complete proof of Lemma 5 by purely combinatorial means. On the other hand, we can give an easy proof of \eqref{N_d} (as well as Lemmas 4 and 5) by appealing to \eqref{ACGHNd}, as follows. 
Namely, the Exponential Formula implies that
\[
\sum_{d \geq 0} N_d(g,m) z^d= \mbox{exp}(\sum_{n >0} [m \phi_1+ (2g-2)\phi_2]z^n)
\]
where $\phi_1$ and $\phi_2$ are rational functions of $n$. It suffices to show that
\beg{equation}\label{phi}
\phi_1= (-1)^{n-1}\fr{\binom{2n-1}{n-1}}{n} \text{ and } \phi_2= (-1)^{n-1} \fr{4^{n-1}- \binom{2n-1}{n-1}}{n}.
\end{equation}

Now let $\widetilde{\pi}=g-1$. Note that \eqref{ACGHNd} implies that
\beg{equation}\label{NewACGHNd}
N_d(g,m)= \sum_{\al=0}^d (-1)^{\al} \binom{\widetilde{\pi}+2d-1-m}{\al}\binom{\widetilde{\pi}+1}{d-\al}.
\end{equation}
We view the expression on the right side of \eqref{NewACGHNd} as a polynomial in $m$ and $\widetilde{\pi}$ with \\coefficients in $\mb{Q}[d]$, whose term of degree 1 in $m$ and degree 0 in $\widetilde{\pi}$ is $\phi_1$, and whose  term of degree 0 in $m$ and degree 1 in $\widetilde{\pi}$ is $\phi_2$. As a matter of notation, 
given any polynomial $Q$ in $m$ and $\widetilde{\pi}$, we let $[m^{\al} \widetilde{\pi}^{\bet}]Q$ denote the coefficient of $m^{\al} \widetilde{\pi}^{\bet}$ in $Q$.

\medskip
To prove the first identity in \eqref{phi}, note that, by \eqref{NewACGHNd},
{\small
\[
\phi_1 = [m] \sum_{\al=0}^n (-1)^{\al} \binom{\widetilde{\pi}+ 2n-1-m}{\al} \binom{\widetilde{\pi}+1}{n-\al}.
\]
}

Similarly, to prove the second identity in \eqref{phi}, note that, by \eqref{NewACGHNd},
{\small
\[
\phi_2 = \fr{1}{2} \cdot [\widetilde{\pi}] \sum_{\al=0}^n (-1)^{\al} \binom{\widetilde{\pi}+ 2n-1-m}{\al} \binom{\widetilde{\pi}+1}{n-\al}.
\]
}

In other words, \eqref{phi} may be reduced to hypergeometric identities, which are handled by the Wilf--Zeilberger algorithm \cite[p.83]{Ko}. To prove the second identity, i.e., that
\[
(-1)^{n-1} \fr{4^{n-1}- \binom{2n-1}{n-1}}{n}= \fr{1}{2} \cdot [\widetilde{\pi}] \sum_{\al=0}^n (-1)^{\al} \binom{\widetilde{\pi}+ 2n-1-m}{\al} \binom{\widetilde{\pi}+1}{n-\al},
\]
one needs to show that
\[
1= \sum_{i=-\infty}^{\infty} \bigg(\fr{n}{i+1}- \fr{n}{i+2}+2 \bigg) \fr{n!i!}{(n+1+i)!}\binom{n-2}{i}.
\]
The interested reader may see \cite{Co1} for details.
\end{proof}

{\fl \bf NB}: As noted above, we have $\fr{2}{(1+4z)^{1/2}+1}=C(-z)$, so \eqref{N_dbis} may be reexpressed in the following more compact form:
\beg{equation}\label{N_dbisbis}
\sum_{d \geq 0} N_d(g,m) z^d= C(-z)^{2g-2-m} \cdot (1+4z)^{\fr{g-1}{2}}.
\end{equation}

\subsection{Generating functions for tautological coefficients when $r=1$}\label{genfunctions}
When $r=1$, the results of the preceding subsections imply that
\beg{equation}\label{explicitP}
\beg{split}
P_{\al} &= \biggl[\fr{m+1-2d}{2g}\biggr] N_d(g,m)- \biggl[\fr{d+1}{2g}\biggr] N_{d+1}(g,m+1), \\
P_{\be} &= \fr{-mP_{\al}+ d N_d(g,m)}{g-1}, \text{ and }\\
P_c &=-N_d(g,m).
\end{split}
\end{equation}

Given our generating function \eqref{N_dbis} for $N_d(g,m)$ (and whence, for $P_c=P_c(d,g,m)$), determining generating functions for $P_{\al}=P_{\al}(d,g,m)$ and $ P_{\be}=P_{\be}(d,g,m)$ is now a purely formal matter. Namely, let
\beg{equation}\label{X,Y}
Z_{g,m}(z):=\biggl(\fr{2}{(1+4z)^{1/2}+1} \biggr)^{2g-2-m} \cdot (1+4z)^{\fr{g-1}{2}}.
\end{equation}
Then, according to \eqref{N_dbis},
\[
\sum_{d \geq 0} N_d(g,m) z^d= Z_{g,m}(z).
\]
By \eqref{explicitP}, it follows that
\beg{equation}\label{Palphagf}
\beg{split}
\sum_{d \geq 0} P_{\al}(d,g,m) z^d&= \fr{1}{2g} \cdot \sum_{d \geq 0} [(m+1-2d) N_d(g,m)- (d+1)N_{d+1}(g,m+1)] z^d \\
&= \biggl(\fr{m+1}{2g}-\fr{z}{g} \fr{d}{dz}\biggr) Z_{g,m}(z)- \biggl(\fr{1}{2g} \cdot \fr{d}{dz} \biggr) Z_{g,m+1}(z) \\
&= Z_{g,m}(z) \biggl[ \fr{1}{2}- \fr{1}{2(1+4z)^{1/2}} \biggr].
\end{split}
\end{equation}
Similarly, we have 
\[
\beg{split}
\sum_{d \geq 0} P_{\be}(d,g,m) z^d &= -\fr{m}{g-1} \sum_{d \geq 0} P_{\al}(d,g,m) z^d + \fr{z}{g-1} \cdot \fr{d}{dz} Z_{g,m}(z)\\
&= Z_{g,m}(z) \biggl[ \fr{2z}{1+4z}- \fr{4z}{(1+4z)^{1/2} ((1+4z)^{1/2}+1)} \biggr].
\end{split}
\]

\subsection{From generating functions to generalized hypergeometric series}\label{hypergeomfns}
Using the results of the preceding subsection, it is possible to realize $P_c$, $P_{\al}$ and $P_{\be}$ as linear combinations of generalized hypergeometric series whenever $r=1$. Namely, we have the following result.

\beg{thm}\label{hypergeom} When $r=1$, the tautological secant-plane divisor coefficients $P_{\al}=P_{\al}(d,g,m)$, $P_{\be}=P_{\be}(d,g,m)$, and $P_c=P_c(d,g,m)$ are given by
\[
P_c=-\fr{g! (2g-2-m)!}{(g-2d)!d!(2g-2-m+d)!} { }_3F_2 \biggl[\beg{array}{ccc} -\fr{g}{2}+ \fr{m}{2}+ 1-d, & -\fr{g}{2}+ \fr{m+3}{2}-d, & -d \\ \fr{g+1}{2}-d, & \fr{g}{2}+1-d \end{array} \biggl| 1\biggr],
\]
\[
P_{\al}=\fr{g! (2g-2-m)!}{2(g-2d)!d!(2g-2-m+d)!} { }_3F_2 \biggl[\beg{array}{ccc} -\fr{g}{2}+ \fr{m}{2}+ 1-d, & -\fr{g}{2}+ \fr{m+3}{2}-d, & -d \\ \fr{g+1}{2}-d, & \fr{g}{2}+1-d \end{array} \biggl| 1\biggr]
\]
\vspace{-3pt}
\[
-\fr{(g-1)!(2g-2-m)!}{2(g-2d-1)!d!(2g-2-m+d)!} { }_3F_2 \biggl[\beg{array}{ccc} -\fr{g}{2}+ \fr{m}{2}+ 1-d, & -\fr{g}{2}+ \fr{m+1}{2}-d, & -d \\ \fr{g+1}{2}-d, & \fr{g}{2}-d \end{array} \biggl| 1\biggr],
\]
\[
P_{\be}= \fr{2(g-2)!(2g-2-m)!}{(g-2d)!(d-1)!(2g-3-m+d)!} { }_3F_2 \biggl[\beg{array}{ccc} -\fr{g}{2}+ \fr{m}{2}+ 1-d, & -\fr{g}{2}+ \fr{m+3}{2}-d, & 1-d \\ \fr{g+1}{2}-d, & \fr{g}{2}+1-d \end{array} \biggl| 1\biggr]
\]
\vspace{-3pt}
\[
-\fr{2(g-1)!(2g-1-m)!}{(g+1-2d)!(d-1)!(2g-2-m+d)!} { }_3F_2 \biggl[\beg{array}{ccc} -\fr{g}{2}+ \fr{m}{2}+ 1-d, & -\fr{g}{2}+ \fr{m+3}{2}-d, & 1-d \\ \fr{g}{2}+1-d, & \fr{g+3}{2}-d \end{array} \biggl| 1\biggr],
\]
\end{thm}

\beg{proof}
Recall (see, e.g., \cite{PWZ}) that 
\[
{}_pF_q \biggl[\beg{array}{ccc} a_1, & \dots,& a_p \\ b_1,& \dots,& b_q \end{array}\bigg| \phi \biggr]= \sum_{k=0}^{\infty} \fr{(a_1)^{(k)} \cdots (a_p)^{(k)}}{(b_1)^{(k)} \cdots (b_q)^{(k)}} \fr{\phi^k}{k!}
\]
where $(u)^{(k)}= \fr{\Ga(u+k)}{\Ga(u)}$ is the Pochhammer symbol. 

Using \eqref{N_dbisbis}, we find that
\[
\beg{split}
P_c(d,g,m)&=-[z^d] (C(-z)^{2g-2-m} \cdot (1+4z)^{\fr{g-1}{2}}) \\
&= -\sum_{k=0}^d [z^k] C(-z)^{2g-2-m} \cdot [z^{d-k}] (1+4z)^{\fr{g-1}{2}}.
\end{split}
\]
Here
\[
[z^k] C(-z)^{2g-2-m}= (-1)^k \fr{2g-2-m}{k+ 2g-2-m} \binom{2k+ 2g-3-m}{k}.
\]
It follows that
\beg{equation}\label{Pcv1}
\beg{split}
P_c(d,g,m) &= -\sum_{k=0}^d (-1)^k \fr{2g-2-m}{k+ 2g-2-m} \binom{2k+ 2g-3-m}{k} 4^{d-k} \binom{\fr{g-1}{2}}{d-k} \\
&= -\sum_{k=0}^{\infty} (-1)^k \fr{2g-2-m}{k+ 2g-2-m} \binom{2k+ 2g-3-m}{k} 4^{d-k} \binom{\fr{g-1}{2}}{d-k}.
\end{split}
\end{equation}
We claim first that the expression on the right side of \eqref{Pcv1} equals
\beg{equation}\label{Pcv2}
-4^d \binom{\fr{g-1}{2}}{d}{ }_3F_2 \biggl[\beg{array}{ccc} g-1-\fr{m}{2}, & g-\fr{m+1}{2}, &-d \\ 2g-1-m, & \fr{g+1}{2}-d \end{array}\biggl| 1\biggr].
\end{equation}
Clearly, the $0$th terms of the two series agree.

Moreover, an easy calculation now shows that for both the right side of \eqref{Pcv1} and the expansion of the hypergeometric series \eqref{Pcv2}, the quotient of the $(k+1)$th over the $k$th term equals
\[
\fr{1}{2} \fr{(2k-1+2g-m)(k-d)(2k+2g-2-m)}{(2k+g+1-2d)(k+1)(k-1+2g-m)}.
\]

On the other hand, for every nonnegative integer $n$, we have the following equality of hypergeometric series \cite{GR}:
\[
{ }_3F_2 \biggl[\beg{array}{ccc} w, & x, & -n \\ y, & z \end{array}\biggl| 1\biggr]= \fr{(-w-x+y+z)_n}{(z)_n} { }_3F_2 \biggl[\beg{array}{ccc} -w+y, & -x+y, &-n \\ y, & -w-x+y+z \end{array}\biggl| 1\biggr].
\]
Taking $w=g-1-\fr{m}{2}$, $x= g-\fr{m+1}{2}$, $y=\fr{g+1}{2}-d$, $z=2g-1-m$, and $n=d$, we deduce that
\[
-4^d \binom{\fr{g-1}{2}}{d}{ }_3F_2 \biggl[\beg{array}{ccc} g-1-\fr{m}{2}, & g-\fr{m+1}{2}, &-d \\ 2g-1-m, & \fr{g+1}{2}-d \end{array}\biggl| 1\biggr]\\
\]
\vspace{-3pt}
\[
= -4^d \binom{\fr{g-1}{2}}{d} \fr{(\fr{g}{2}+1-d)_d}{(2g-1-m)_d} { }_3F_2 \biggl[\beg{array}{ccc} -\fr{g}{2}+ \fr{m}{2}+ 1-d, & -\fr{g}{2}+ \fr{m+3}{2}-d, & -d \\ \fr{g+1}{2}-d, & \fr{g}{2}+1-d \end{array} \biggl| 1\biggr]
\]
Finally, it is elementary to check that
\[
4^d \binom{\fr{g-1}{2}}{d} \fr{(\fr{g}{2}+1-d)_d}{(2g-1-m)_d}= \fr{g! (2g-2-m)!}{(g-2d)!d!(2g-2-m+d)!};
\]
it follows that
\[
P_c=-\fr{g! (2g-2-m)!}{(g-2d)!d!(2g-2-m+d)!} { }_3F_2 \biggl[\beg{array}{ccc} -\fr{g}{2}+ \fr{m}{2}+ 1-d, & -\fr{g}{2}+ \fr{m+3}{2}-d, & -d \\ \fr{g+1}{2}-d, & \fr{g}{2}+1-d \end{array} \biggl| 1\biggr],
\]
as desired.
The proofs of the the other equalities are similar.
\end{proof}

\subsection{The case $r=s$} From \cite[Ch. VIII, Prop. 4.2]{ACGH}, we see that the expected number of $2s$-secant $(s-1)$-planes to a curve $C$ of degree $(m+1)$ and genus $g$ in $\mb{P}^{s+1}$ is
\[
A^{\pr}=\fr{(-1)^{\binom{s}{2}}}{2} [((1+t_1)(1+t_2))^{m-g-s}(1+t_1+t_2)^g (t_1-t_2)^2]_{t_1^{s+1}t_2^{s+1}}.
\]

Similarly, Macdonald's formula \cite[Thm. 4]{M} specializes nicely in the case $r=s$. It implies that the expected number of $(2s-1)$-secant $(s-1)$-planes to a curve $C$ of degree $m$ and genus $g$ in $\mb{P}^{s+1}$ that intersect a disjoint line is
\[
A=\fr{(-1)^{\binom{s}{2}}}{2} [((1+t_1)(1+t_2))^{m-g-s}(1+t_1+t_2)^g (t_1-t_2)^2(2t_1t_2+t_1+t_2)]_{t_1^{s+1} t_2^{s+1}}.
\]
To see this, simply note that the condition imposed by requiring an $(s-1)$-plane to intersect a line in $\mb{P}^{s+1}$ defines the Schubert cycle $\sig_1$ in $\mb{G}(s-1,s+1)$; the formula for $A$ above follows from Macdonald's by a straightforward calculation.

It is natural to wonder whether the numbers $A$ and $A^{\pr}$ admit neat combinatorial descriptions when $r=s$ (or indeed, whenever $1< r \leq s$, just as they do when $r=1$.
 
\section{Planes incident to linear series on a general curve when $\rho=1$}\label{rho=1}
In this section, we use the results of the preceding one to deduce a new formula for the number $N^{\pr, d-r-1}_d$ of linear series with exceptional secant planes on a general curve of genus $g$, which is applicable whenever $\rho=1$ and $\mu=-1$. (By Theorem~\ref{thm3}, that number is always finite.) Namely, we have the following result.
\beg{thm}\label{genericcurve}
Let $\rho=1, \mu=-1$. The number $N^{\pr, d-r-1}_d$ of linear series $g^s_m$ with $d$-secant $(d-r-1)$-planes on a general curve of genus $g$ is given by
\vspace{-3pt}
\[
\beg{split}
&N^{\prime, d-r-1}_d= \fr{(g-1)! 1! \cdots s!}{(g-m+s)! \cdots (g-m+2s-1)!(g-m+2s+1)!} \\
&[(-gm+2gs+m^2-3ms+2s^2-m+s+g)A+ (gd+g-md-m+2sd+2s+d+1)A^{\pr}].
\end{split}
\]
where $A$ and $A^{\pr}$ are as defined in Section~\ref{testfamilies}.
\end{thm}

\beg{proof}
We use the basic set-up of Section~\ref{testfamilies}, as well as the relations among the tautological coefficients $P_{\al}, P_{\be}, \text{ and } P_{\de_0}$ obtained there, to prove Theorem~\ref{genericcurve}. Namely, let $C$ denote a general curve of genus $g$ such that $\rho(g,s,m)=1$, and consider the test family $\pi:\mc{X} \ra B$ with total space $\mc{X}=W^s_m(C) \times C$ and base $B=W^s_m(C)$. Let $\mc{L}$ denote the pullback of any degree-$m$ Poincar\'{e} bundle $\widetilde{\mc{L}} \ra \mbox{Pic}^m(C) \times C$ by
  the inclusion $i \times 1_C: W^s_m(C) \times C \ra \mbox{Pic}^m(C) \times C$. Let $\theta$ and $\eta$ denote the integral cohomology classes of the pullbacks to $\mbox{Pic}^m(C) \times C$ of the theta divisor on $\mbox{Pic}^m(C)$ and a point on $C$, respectively. As explained in \cite[Ch. $8$]{ACGH}, we then have
\[
\beg{split}
c_1(\mc{L}) &= (m \eta+ \ga) \cdot \nu^* w^s_m \\
&= \ub{\De(g-m+s,\dots,g-m+s)}_{(s+1) \text{ times}} \cdot (m \eta+ \ga) \cdot \nu^* \theta^{g-1}
\end{split}
\]
where $\nu: \mbox{Pic}^m(C) \times C \ra \mbox{Pic}^m(C)$ is the natural projection, and where
$\De(a_1,\dots,a_n)$ denotes the determinant of the $n \times n$ matrix with $(i,j)$th entry $\fr{1}{(a_i+j-i)!}$ (in our case, $n=s+1$).

It follows immediately that
\[
\beg{split}
\al&= \De(g-m+s,\dots,g-m+s) \cdot (m \eta+ \ga)^2 \cdot \nu^*
\theta^{g-1} \\
&= \De(g-m+s,\dots,g-m+s) \cdot (-2\eta \theta) \cdot \nu^*
\theta^{g-1} \\
&= -2 g! \De(g-m+s,\dots,g-m+s)
\end{split}
\]
and, likewise, that
\[
\be = \ga= \de=0,
\]
since $\om= \pi_2^* K_C= \pi_2^* (2g-2) \{\mbox{pt}_C\}$ in this case.

Finally, let  $\Ga$ denote any section of $\mbox{Pic}^d(C) \times C \ra \mbox{Pic}^d(C)$ associated to a divisor of large degree on $C$. Note that $\mc{V}$ is the kernel bundle for the evaluation map
\[
\mc{E}:=\nu_*(\widetilde{\mc{L}}(\Gamma)) \st{\mbox{ev}}{\ra} \nu_*(\widetilde{\mc{L}}(\Gamma)/\widetilde{\mc{L}})=:\mc{F}
\]
of vector bundles over $\mbox{Pic}^d(C)$, restricted to the
locus along which $\mbox{ev}$ has a kernel of rank $(s+1)$. On the other hand, the vector bundle $\mc{F}$ has trivial Chern classes. Accordingly,
the kernel number formula of \cite{HT} yields
\[
\beg{split}
c &= - \Delta_{g-m+s+1,g-m+s, \dots, g-m+s}(c_t(-\mc{E})) \\
&= -g! \Delta(g-m+s+1,g-m+s, \dots, g-m+s).
\end{split}
\]
Here $\De(a_1,\dots,a_n)(\mc{F})$ denotes the determinant of the $n \times n$ matrix with $(i,j)$th entry $c_{a_i+j-i}(\mc{F})$, for any vector bundle $\mc{F}$.

On the other hand, from the results of Section~\ref{testfamilies}, we see that
\[
P_{\al}= \fr{(m-s)A-(d+1)A^{\pr}}{2g} \text{ and } P_c=-A.
\]
It follows immediately that
\[
\beg{split}
N^{\prime, d-r-1}_d&= -2 g! \De(g-m+s,\dots,g-m+s) \fr{(m-s)A-(d+1)A^{\pr}}{2g} \\
&+ g! \Delta(g-m+s+1,g-m+s, \dots, g-m+s) A.
\end{split}
\]
To simplify the latter expression, we use the well-known fact (see, e.g., \cite[p.320]{ACGH}) that
\[
\De(a_1,\dots,a_n)= \fr{\prod_{j>i} (a_i-a_j-i+j)}{\prod_{i=1}^n (a_i-i+n-1)!}.
\]
We deduce that
\[
\beg{split}
&\De(g-m+s,\dots,g-m+s)= \fr{s! \cdots 1!}{(g-m+2s)! \cdots (g-m+s)!} \text{ and }\\
&\De(g-m+s+1,g-m+s,\dots,g-m+s) \\
&= \fr{(s+1)!(s-1)! \cdots 1!}{(g-m+2s+1)!(g-m+2s-1)! \cdots (g-m+s)!}.
\end{split}
\]
It follows that
\[
\beg{split}
&N^{\prime, d-r-1}_d= \fr{(g-1)! 1! \cdots s!}{(g-m+s)! \cdots (g-m+2s-1)!(g-m+2s+1)!} \\
&\biggl[-2g(g-m+2s+1) \biggl(\fr{(m-s)A}{2g}- \fr{(d+1)A^{\pr}}{2g} \biggr)+ g(s+1)A \biggr] \\
&= \fr{(g-1)! 1! \cdots s!}{(g-m+s)! \cdots (g-m+2s-1)!(g-m+2s+1)!} \\
&[(-gm+2gs+m^2-3ms+2s^2-m+s+g)A+ (gd+g-md-m+2sd+2s+d+1)A^{\pr}].
\end{split}
\]
\end{proof}

\subsection{The case $r=1$}
Following our usual practice, we now specialize to the case $r=1$, so that $N^{\prime, d-r-1}_d$ counts $(2d-1)$-dimensional series with $d$-secant $(d-2)$-planes. 
Here we obtain stronger results by applying Theorem~\ref{hypergeom}, which characterizes the tautological secant-plane coefficients $P=P(d,g,m)$ in terms of hypergeometric series. Because $\rho=1$, we have
\[
g=2ad+1, \text{ and }m= (a+1)(2d-1)+1
\]
for suitably chosen positive integers $a$ and $d$ (here, as usual, $d$ denotes incidence).

Accordingly, we have
{\small
\beg{equation}\label{ndprimead}
\beg{split}
N^{\prime, d-2}_d &=\fr{g! 1! \cdots s!}{(g-m+s)! \cdots (g-m+2s-1)!(g-m+2s+1)!} [-(s+1)P_c- 2(g-m+2s+1)P_{\al}] \\
&= \fr{(2ad+1)! 1! \cdots (2d-1)!}{a! \cdots (a+2d-2)!(a+2d)!} [-2d P_c(a,d) -2(2d+a) P_{\al}(a,d)] \\
&=\fr{(2ad+1)! 1! \cdots (2d-1)!}{a! \cdots (a+2d-2)!(a+2d)!} \bigg[-2d P_c(a,d)- 2(2d+a) \bigg(-\fr{1}{2}P_c(a,d)+ P_{\al,2}(a,d)\bigg)\bigg] \text{ by Theorem 5}\\
&=\fr{(2ad+1)! 1! \cdots (2d-1)!}{a! \cdots (a+2d-2)!(a+2d)!} [aP_c(a,d)- (4d+2a)P_{\al,2}(a,d)]
\end{split}
\end{equation}
}
where
{\small
\[
P_c(a,d)= -\fr{(2ad+1)![(2a-2)d+a]!}{[(2a-2)d+1]!(2ad-d+a)!d!}  { }_3F_2 \bigg[\beg{array}{ccc} \fr{1}{2}- \fr{a}{2}, & 1-\fr{a}{2}, & -d \\ (a-1)d+1, & (a-1)d+ \fr{3}{2} \end{array} \bigg| 1\bigg],
\]
}
i.e.,
{\small
\beg{equation}\label{Pc}
P_c(a,d)= -\fr{(2ad+1)!}{(2ad-d+a)!d!} \sum_{i=0}^{\lfl \fr{a-1}{2} \rfl} (-1)^i \fr{((2a-2)d+a)!}{((2a-2)d+2i+1)!} \cdot \fr{d!}{(d-i)!} \cdot \fr{(a-1)!}{(a-1-2i)!} \cdot \fr{1}{i!},
\end{equation}
and
{\small
\[
P_{\al,2}(a,d)= \fr{(2ad+1)![(2a-2)d+a-1]!}{[(2a-2)d-1]!(2ad-d+a)!(d+1)!} { }_3F_2 \bigg[\beg{array}{ccc} - \fr{a}{2}, & \fr{1}{2}-\fr{a}{2}, & -d \\ (a-1)d, & (a-1)d+ \fr{1}{2} \end{array} \bigg| 1\bigg],
\]
}
i.e.,
{\small}
\beg{equation}\label{Palpha}
P_{\al,2}(a,d)=-\fr{(2ad)!}{2(2ad-d+a)!d!} \sum_{i=0}^{\lfl \fr{a}{2} \rfl} (-1)^i \fr{((2a-2)d+a)!}{((2a-2)d+2i)!} \cdot \fr{d!}{(d-i)!} \cdot \fr{a!}{(a-2i)!} \cdot \fr{1}{i!}.
\end{equation}
}
As a consequence, we deduce the following result, which establishes that the secant-plane loci in question are nearly always non-empty. Note that this is far from obvious from Theorem~\ref{genericcurve}, as the formulas for $A$ and $A^{\pr}$ given in \cite{ACGH} involve alternating sums of binomial coefficients.
\beg{thm}\label{ndprimepos}
The number of series with exceptional secant planes $N^{\prime, d-2}_d$ is zero when either $a=1$ or $d=1$, and is positive whenever $a>1$ and $d>1$.
\end{thm}

\beg{proof}
First assume that $a=1$. Note that \eqref{Pc} implies that
\[
P_c(1,d)= -\fr{(2d+1)!(d+1)}{((d+1)!)^2}, \text{ and } P_{\al,2}(1,d)= -\fr{(2d)!(d+1)}{2((d+1)!)^2}.
\]
It follows that 
\[
P_c(1,d)- (4d+2)P_{\al,2}(1,d)=0;
\]
and whence, by \eqref{ndprimead}, that $N^{\prime, d-2}_d=0$.

Similarly, \eqref{Palpha} implies that
\[
P_c(a,1)=-(a+2), \text{ and } P_{\al,2}(a,1)= -\fr{a}{2},
\]
so that
\[
P_c(a,1)- (4d+2)P_{\al,2}(a,1)=0,
\]
and, therefore, $N^{\prime, d-2}_d=0$.

Now assume that $a>1$, and $d>1$. In view of \eqref{ndprimead}, we need only show that $(4d+2a) P_{\al,2}(a,d)< aP_c(a,d)$ whenever $a>1$ and $d>1$, i.e., that
\[
\beg{split}
&\sum_{i=0}^{\lfl \fr{a-1}{2} \rfl} (-1)^i \fr{((2a-2)d+a)!}{((2a-2)d+2i+1)!} \cdot \fr{d!}{(d-i)!} \cdot \fr{a!}{(a-1-2i)!} \cdot \fr{2ad+1}{i!} \\
&< \sum_{i=0}^{\lfl \fr{a}{2} \rfl} (-1)^i \fr{((2a-2)d+a)!}{((2a-2)d+2i)!} \cdot \fr{d!}{(d-i)!} \cdot \fr{a!}{(a-2i)!} \cdot \fr{2d+a}{i!}.
\end{split}
\]
To this end, write
\[
\beg{split}
P_1^{(i)}&= \fr{((2a-2)d+a)!}{((2a-2)d+2i)!} \cdot \fr{d!}{(d-i)!} \cdot \fr{a!}{(a-2i)!} \cdot \fr{2d+a}{i!} \text{ for all } 0 \leq i \leq \lfl \fr{a}{2} \rfl, \text{ and }\\
P_2^{(i)}&= \fr{((2a-2)d+a)!}{((2a-2)d+2i+1)!} \cdot \fr{d!}{(d-i)!} \cdot \fr{a!}{(a-1-2i)!} \cdot \fr{2ad+1}{i!} \text{ for all }0 \leq i \leq \lfl \fr{a-1}{2} \rfl.
\end{split}
\]
Here
\[
\fr{P_1^{(i)}}{P_2^{(i)}}= \fr{(2d+a)((2a-2)d+2i+1)}{(a-2i)(2ad+1)}
\]
for all $i \leq \lfl \fr{a-1}{2} \rfl$. So unless $a=4k+2$ for some $k \geq 1$ (a case we will handle separately below), we need only show that the quantity
\[
T:=\sum_{i=0}^{\lfl \fr{a-1}{2} \rfl} (-1)^i \biggl[ \fr{(2d+a)((2a-2)d+2i+1)}{(a-2i)(2ad+1)}-1 \biggr] P_2^{(i)}
\]
is positive. (When $a=4k+2$, and only in that case, we have 
\[
\sum_{i=0}^{\lfl \fr{a}{2} \rfl} (-1)^i P_1^{(i)}- \sum_{i=0}^{\lfl \fr{a-1}{2} \rfl} (-1)^i P_2^{(i)}<T.)
\]
To this end, it suffices, in turn, to show that $\fr{T_i}{T_{i+1}}>1$, where
\[
T_i= \biggl[ \fr{(2d+a)((2a-2)d+2i+1)}{(a-2i)(2ad+1)}-1 \biggr] P_2^{(i)}.
\]
Simplifying $\fr{T_i}{T_{i+1}}$ yields
\[
\fr{T_i}{T_{i+1}}= \fr{2[(2a-2)d^2+(2ai+2i-a+1)d+ (a+1)i][(a-1)d+i+1][(2a-2)d+2i+3](i+1)}{[(2a-2)d^2+(2ai+2i+a+3)d+(a+1)(i+1)](a-1-2i)(a-2i)(d-i)}.
\]
Here 
\[
\beg{split}
&\fr{(2a-2)d^2+(2ai+2i-a+1)d+ (a+1)i}{(2a-2)d^2+(2ai+2i+a+3)d+(a+1)(i+1)} \\
&= 1-\fr{(2a-2)d+(a+1)}{(2a-2)d^2+(2ai+2i+a+3)d+(a+1)(i+1)}> 1- \fr{1}{d},
\end{split}
\]
while
\[
\fr{((2a-2)d+2i+3}{a-1-2i}>2, \text{ and } \fr{(a-1)d+i+1}{(a-2i)(d-i)}> 1-\fr{1}{a}.
\]
We conclude that $\fr{T_i}{T_{i+1}}>1$ whenever $a>1$ and $d>1$. 

It remains to treat the case where $a=4k+2$ for some $k \geq 1$. To conclude the proof of our theorem, it will suffice to show that
\beg{equation}\label{4kdiff}
P_1^{(i)}- P_1^{(i+1)}- P_2^{(i)}>0
\end{equation}
for $i=\fr{a}{2}-2$.
Simplifying yields
\[
P_1^{(i)}- P_1^{(i+1)}- P_2^{(i)}= \fr{((2a-2)d+a)! d! a!}{(a-2i)!(d-i)!(i+1)!((2a-2)d+2i+2)!} \widetilde{Q}(a,d,i)
\]
where
\[
\beg{split}
\widetilde{Q}(a,d,i) &=(8a^2i-16ai+ 8a^2-16a+ 8i+8) d^3 \\
&+(-12+8a^2i^2+4a^2i+32ai+8ai^2+18a-32i-6a^2-24i^2)d^2 \\
&+(-4a+4+8ai^3-a^3+10a^2i-4ai+4a^2i^2+12i+a^2+20i^2+16i^3) d \\
&+(4i+a^3i-4a^2i^2+4ai+10ai^2+8ai^3-a^2i+8i^2+4i^3).
\end{split}
\]
Taking $i=\fr{a}{2}-2$, we find that
\[
\beg{split}
\widetilde{Q}\biggl(a,d,\fr{a}{2}-2\biggr) &=(4a^3-16a^2+20a-8) d^3+ (2a^4-12a^3+12a^2+18a-44) d^2\\
&+ (2a^4-14a^3+ 24a^2+ 2a-68)d + \biggl(\fr{1}{2}a^4- \fr{7}{2}a^3+ 12a^2-22a-8\biggr),
\end{split}
\]
which is positive whenever $a \geq 4$ and $d \geq 2$.
\end{proof}
Finally, we calculate the asymptotic behavior of $N^{\prime, d-2}_d$, using \eqref{ndprimead}. To this end, note that \eqref{ndprimead} implies that when $r=1$,
{\small
\[
\beg{split}
P_c(a,d)&= -\fr{(2ad+1)!}{(2ad-d+a)!d!} \biggl[ \fr{((2a-2)d+a)!}{((2a-2)d+1)!}+ O(d^{a-2}) \biggr], \text{ and } \\
P_{\al,2}(a,d)&=-\fr{(2ad)!}{2(2ad-d+a)!d!} \biggl[ \fr{((2a-2)d+a)!}{((2a-2)d)!}+ O(d^{a-1}) \biggr].
\end{split}
\]
}
It follows that
\[
\beg{split}
N^{\prime, d-2}_d &=\fr{(2ad+1)! 1! \cdots (2d-1)!}{a! \cdots (a+2d-2)!(a+2d)!} [aP_c(a,d)- (4d+2a)P_{\al,2}(a,d)] \\
&=\fr{(2ad+1)! 1! \cdots (2d-1)!}{a! \cdots (a+2d-2)!(a+2d)!} \cdot \fr{(2ad)!}{(2ad-d+a)!d!} \biggl[ \fr{2d((2a-2)d+a)!}{((2a-2)d)!} \\
&-\biggl(\fr{(2ad+1)a ((2a-2)d+a)!}{((2a-2)d+1)!}+ \fr{a((2a-2)d+a)!}{((2a-2)d)!} \biggr) + O(d^{a-1}) \biggr] \\
&=\fr{(2ad+1)! 1! \cdots (2d-1)!}{a! \cdots (a+2d-2)!(a+2d)!} \cdot \fr{(2ad)!}{(2ad-d+a)!d!} \cdot \fr{((2a-2)d+a)!}{((2a-2)d+1)!} \\
&[(4a-4)d^2+ (-4a^2+2a+2)d+O(1)].
\end{split}
\]

{\fl \bf NB:} 
Theorem~\ref{ndprimepos} establishes that no series with $a=1$ and $\rho=1$ on a general curve $C$ of genus $g$ admits $d$-secant $(d-2)$-planes. This is easy to explain on geometric grounds. Namely, every such series $g^{2d-1}_{4d-1}$ can be realized as a subseries of a canonical series $g^{2d}_{4d}$ with a base point, so its image necessarily arises as the image of a canonical curve $\widetilde{C} \sub \mb{P}^{2d}$ under projection from a point along $\widetilde{C}$. Moreover, $d$-secant $(d-2)$-planes of our original series are in bijection with $(d+1)$-secant $(d-1)$-planes to $\widetilde{C}$. On the other hand, any $(d+1)$-secant $(d-1)$-plane to $\widetilde{C}$ defines an inclusion of linear series
\beg{equation}\label{canonicalinclusion}
g^d_{3d-1}+p_1+\dots+p_{d+1} \hra g^{2d}_{4d}
\end{equation}
along $C$. But in fact $\rho(2d+1,d,3d-1)<0$; whence, by the Brill--Noether theorem, no inclusions \eqref{canonicalinclusion} exist.

Similarly, when $d=1$, $N^{\prime, d-2}_d$ counts one-dimensional series with base points. Theorem~\ref{ndprimepos} establishes that no such series exist on a general curve of genus $g$, which also confirms the Brill--Noether theorem in a special case.

{\fl \bf \small Dept of Mathematics and Statistics, Queen's University, Kingston, ON K7L 3N6 \\{\it Email address}: \url{cotteril@mast.queensu.ca}}
\end{document}